\documentclass[11pt,oneside]{amsart}
\usepackage{amsmath,amsfonts,amsgen,amstext,amsbsy,amsopn,amsthm, amssymb}
\usepackage{multirow, verbatim, cancel, mathabx}
\usepackage{bbold}
\usepackage{enumitem}
\usepackage[T1]{fontenc}
\usepackage[colorlinks=true,hyperindex=true]{hyperref}
\usepackage{tikz}
\usetikzlibrary{calc}
\usepackage{verbatim}
\usetikzlibrary{arrows,chains,matrix,positioning,scopes}
\usepackage{mathrsfs}
\usepackage{color,soul}
\usepackage{xcolor}
\usepackage{tabulary}
\usepackage{pgfplots}
\usepackage{mathrsfs}
\usepackage{bbm}
\usepackage{dsfont}
\usepackage{bm}
\usepackage{mathtools}
\usepackage{float}

\usepackage{arydshln}
\usepackage{color,soul}
\usepackage{subcaption}
 \usepackage{soul,xcolor}
\usepackage{wrapfig}
\usepackage{lscape}
\usepackage{rotating}
\usepackage{epstopdf}
\usepackage{float}
\usepackage{flafter}

\usepackage{array}
\newcolumntype{R}[1]{>{\raggedleft\arraybackslash}p{#1}}

\usepackage{tabu} 
\usepackage{booktabs}

\usepackage{mathtools}

\usepackage[mathscr]{euscript}

\DeclareSymbolFont{rsfs}{U}{rsfs}{m}{n}
\DeclareSymbolFontAlphabet{\mathscrsfs}{rsfs}

\makeatletter

\usepackage{tkz-fct} 
\usepackage{xcolor}

\usepackage{tkz-fct} 
\usepackage{pgfplots}
\usepackage{setspace}

\usetikzlibrary{matrix}
\usepackage{tikz}

\numberwithin{equation}{section}

\newcounter{smallarabics}

\newcounter{smallroman}
\newenvironment{romanenumerate}
{\begin{list}{{\normalfont\textrm{(\roman{smallroman})}}}
  {\usecounter{smallroman}\setlength{\itemindent}{0cm}
   \setlength{\leftmargin}{5ex}\setlength{\labelwidth}{4ex}
   \setlength{\topsep}{0.75\parsep}\setlength{\partopsep}{0ex}
   \setlength{\itemsep}{0ex}}}
{\end{list}}

\newcommand{\ben}{\begin{romanenumerate}}  
\newcommand{\een}{\end{romanenumerate}}  

\newtheorem{theoreme}{theorem }[section]

\newtheorem{conjecture}[theoreme]{Conjecture}

\newcolumntype{L}{>{\centering\arraybackslash}m{3cm}}

\newcommand\nn\nonumber
\renewcommand\leq\varleq
\renewcommand\geq\vargeq

 \newcommand{\N}{\mathbb{N}}

\newcommand{\E}{\mathcal{E}}

\renewcommand{\epsilon}{\varepsilon}

\usepackage{tikz-cd}

\usepackage{xcolor}
\usepackage{dcolumn}
\usepackage{tabu}

\newcolumntype{A}{D{.}{.}{2.3}}

\usepackage{tikz,tkz-tab}
\usepackage{pgfplots}
\pgfplotsset{compat=1.11}
\usetikzlibrary{positioning}
\makeatletter
      \def\@setcopyright{}
      \def\serieslogo@{}
      \makeatother
\begin{document}

\author{Sylvain Gol\'enia and Marc-Adrien Mandich}
   \address{Univ. Bordeaux, CNRS, Bordeaux INP, IMB, UMR 5251,  F-33400, Talence, France}
   \email{sylvain.golenia@math.u-bordeaux.fr}
      \address{Independent researcher, Jersey City, 07305, NJ, USA}
	\email{marcadrien.mandich@gmail.com}
   

   \title[LAP for discrete Schr\"odinger operator]{Additional numerical and graphical evidence to support some Conjectures on discrete Schr{\"o}dinger operators with a more general long range condition}

   \begin{abstract}
  This document contains additional numerical and graphical evidence to support some of the conjectures mentioned in \cite{GM3}. We give more evidence for $\kappa=3,4$ in dimension 2. As mentioned in that article we still don't quite understand the sets $\boldsymbol{\mu}_{\kappa}(\Delta)$ and $\boldsymbol{\Theta}_{\kappa}(\Delta)$ on $(0,1/2)$ for $\kappa=3$ in dimension 2. Here we give a bunch more threshold energies for $\kappa = 3$ in dimension 2.   
  \end{abstract}

%

 

\maketitle
\hypersetup{linkbordercolor=black}
\hypersetup{linkcolor=blue}
\hypersetup{citecolor=blue}
\hypersetup{urlcolor=blue}
\tableofcontents

\section{Introduction}

The notation in this document is 100\% consistent with that of \cite{GM3}. We start by listing the 2 Conjectures of \cite{GM3} for which the material presented in this document is relevant. The first has to do with the rate of convergence of the thresholds $\E_n \in J_2 = J_2 (\kappa) := \left( 2 \cos(\pi / \kappa), 1+\cos(\pi / \kappa) \right)$.

\begin{conjecture}
\label{conjecture12}
Let $\{ \E_n \}$ be the sequence in \cite[Theorem 1.7]{GM3}. Then $\E_n - \inf J_2 = c(\kappa)/n^2 + o(1/n^2)$, $\forall \kappa \geq 3$, where $c(\kappa)$ means a constant depending on $\kappa$.
\end{conjecture}

For this conjecture, see section \ref{RateConv}. The second conjecture has to do with the existence of a conjugate operator giving a strict Mourre estimate on bands $(\E_n, \E_{n-1}) \subset J_2$.

\begin{conjecture}
\label{conjecture22}
Fix $\kappa \geq 2$. Let $\{ \E_n \}$ be the sequence in \cite[Theorem 1.7]{GM3}. For each interval $( \mathcal{E}_n , \mathcal{E}_{n-1} )$, $n \geq 1$, $\exists$ a conjugate operator $\mathbb{A}(n) = \sum_{q=1} ^{N(n)} \rho_{j_q \kappa} (n) A_{j_q \kappa}$, $A_{j_q \kappa} = \sum_{1 \leq i \leq 2} A_i (j_q,\kappa)$, such that the Mourre estimate for $\Delta$ holds wrt.\  $\mathbb{A}(n)$, $\forall E \in ( \mathcal{E}_n , \mathcal{E}_{n-1} )$. $\mathbb{A}(n)$ is typically not unique. It can be chosen so that $N(n) = 2n$. In particular, $\{\E_n\} = J_2 \cap \boldsymbol{\Theta}_{\kappa}(\Delta)$, $\forall \kappa \geq 2$.
\end{conjecture}

For this conjecture, see sections \ref{sectionk2}, \ref{sectionk3}, \ref{sectionk4} and \ref{Appn=1}. Section \ref{k4m1} lists some thresholds $\in \boldsymbol{\Theta}_{m=1,\kappa=4}(\Delta)$ in dimension 2, and section \ref{section work in progress} gives a graphical illustration of some thresholds $\in \boldsymbol{\Theta}_{m,\kappa=3}(\Delta) \cap (0,1/2)$ in dimension 2, $m \geq 1$.

\section{More evidence for the conjecture on the rate of convergence of $\E_n$}
\label{RateConv}

\begin{figure}[H]
  \centering
 \includegraphics[scale=0.265]{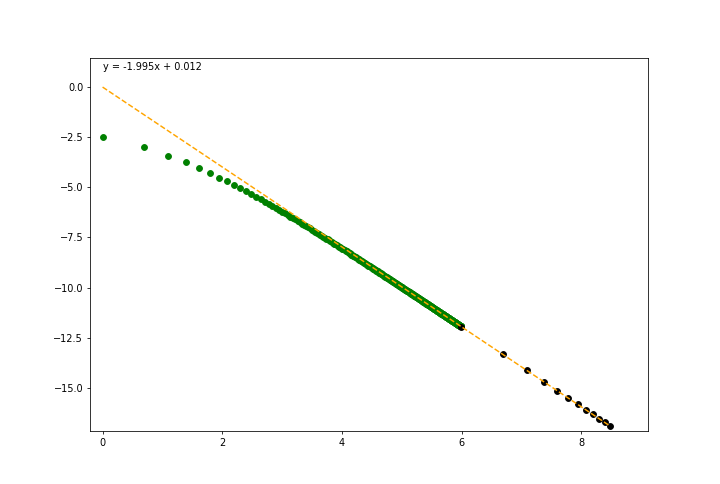}
   \includegraphics[scale=0.22]{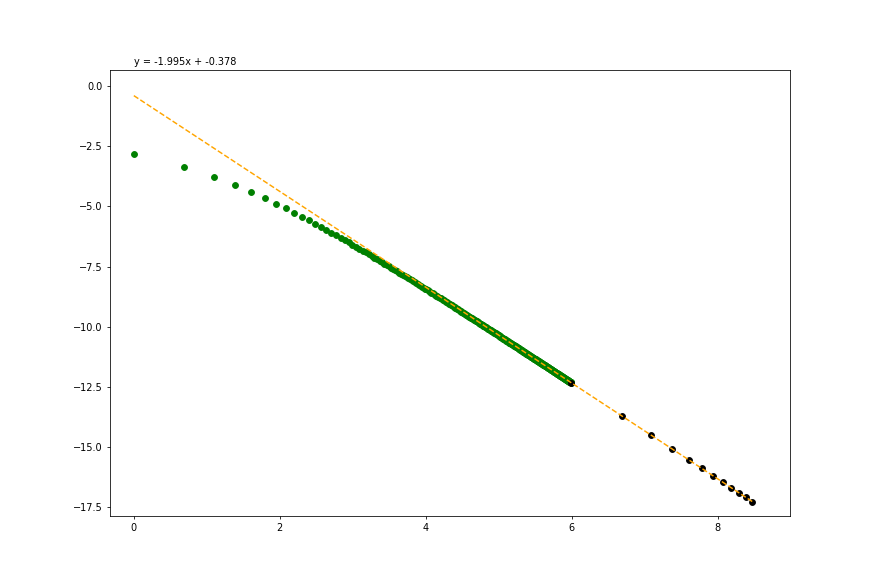}
     \includegraphics[scale=0.27]{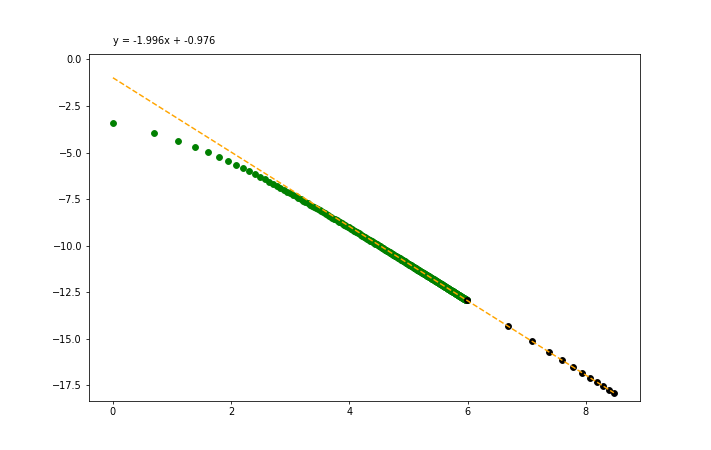}
\caption{Graphs with $\log(n)$ on $x$-axis and $\log(\E_{2n}(\kappa) - 2\cos(\pi / \kappa))$ on $y$-axis. Left : $\kappa=5$ ; Middle : $\kappa=6$, Right : $\kappa=8$. Green dots are $1 \leq n \leq 400$ ; black dots are $n= 400, 800, 1200,...,4800$. Orange line is trend line based on linear regression of black dots.}
\label{fig:test_T3k3_converge}
\end{figure}

Figure \ref{fig:test_T3k3_converge} illustrates solutions $\E_{2n}(\kappa)$ of \cite[Proposition 7.2]{GM3} for $\kappa = 5,6,8$. In these graphs the slope of the orange trend line is close to $-2$, and this is our rationale behind \cite[Conjecture 1.9]{GM3}. 

\section{Energy solutions in $J_2(\kappa) := (2\cos(\pi / \kappa), 1+\cos(\pi / \kappa))$}

Table \ref{table with endpoints3and4} below lists the first few energy solutions for $\kappa \in \{2,3,4\}$ in \cite[Theorem 1.7]{GM3}. In particular we have some exact expressions for $\kappa=4$.
\begin{table}[H]
\small
  \begin{center}
    \begin{tabular}{c|c|c|c|c|c|c|c} 
    $\kappa$ & $\E_5$ & $\E_4$ & $\E_3$ & $\E_2$ & $\E_1$ & $\E_0$   \\ [0.1em]
      \hline
2 & 2/7 & 1/3 & 2/5 & 1/2 & 2/3 & 1 \\ [0.2em]    
3 & $\simeq 1.112$ & $\simeq 1.137$ & $\simeq 1.173$ & $\frac{9+\sqrt{33}}{12} \simeq 1.228$  & $\frac{5+3\sqrt{2}}{7} \simeq 1.320$ & 3/2 \\ [0.2em]
4 & \eqref{E_L_k4_5}  $\simeq 1.476$ & $\simeq 1.491$ & \eqref{E_L_k4_3} $\simeq 1.512$ & \eqref{E_L_k4_2} $\simeq 1.545$ & 8/5 & $1+1/ \sqrt{2} \simeq 1.707$
    \end{tabular}
  \end{center}
    \caption{First values of $\{\E_n\}$ in \cite[Theorem 1.7]{GM3}. $\Delta$ in dimension $2$. $\kappa = 2,3,4$.}
        \label{table with endpoints3and4}
\end{table}
\normalsize


\section{Thresholds in dimension 2, for $\kappa=4$, $m=1$ (the brute force approach)}
\label{k4m1}
In \cite[subsection 4.2]{GM3}, we listed Ansatzes (1) -- (6) to find threshold solutions $\in \boldsymbol{\Theta}_{1,\kappa} (\Delta)$, and in \cite[Lemma 4.2]{GM3} we listed those solutions for $\kappa=3$. Table \ref{sol_kappa_2346} below lists the corresponding solutions for $\kappa=4$. We recall that these solutions belong to $\boldsymbol{\Theta}_{1,\kappa=4} (\Delta)$, but it is an open question for us if they constitute all of $\boldsymbol{\Theta}_{1,\kappa=4} (\Delta)$.

\begin{table}[H]
\small
  \begin{center}
    \begin{tabular}{c|c|c} 
    $\kappa$ & $4$     \\ [0.5em]
      \hline
(1) &  $E=1/2+1/(2\sqrt{2}) \simeq 0.853$, $Y_1=\sqrt{2}/2$  \\ [0.5em]
      & $E= \frac{2 + \sqrt{2} + \sqrt{2 + 4 \sqrt{2}}}{4} \simeq 1.545$, $Y_1 = \sqrt{2}/2$  \\ [0.5em]
    &  $E=\frac{2 + \sqrt{2} - \sqrt{2 + 4 \sqrt{2}}}{4} \simeq 0.161$, $Y_1=\sqrt{2}/2$  \\ [0.5em]
    &  $E=1/2-1/(2\sqrt{2}) \simeq 0.146$, $Y_1=-\sqrt{2}/2$  \\ [0.5em]
 & $E=8/5=1.6$, $Y_1 = E/2 = 0.8$  \\ [0.5em]
(2) &   \\ [0.5em]
(3) & $E=\sqrt{2/5} \simeq 0.632$, $Y_1 = 3/ \sqrt{10} \simeq 0.948$   \\ [0.5em]
&  $E=\sqrt{2}/3 \simeq 0.471$, $Y_1=\sqrt{2}/2$  \\ [0.5em]
(4) &  $E=2/\sqrt{5} \simeq 0.894$, $Y_1 = 0$, $Y_0 = E/2 \simeq 0.447$ \\ [0.5em]
   & $E=1/\sqrt{5} \simeq 0.447$, $Y_1 = 0$, $Y_0=2E \simeq 0.894$  \\ [0.5em]
  & $E = (\sqrt{3}-1)/(2\sqrt{2}) \simeq 0.258$, $Y_1 =0$, $Y_0 = -\sqrt{2}/2$  \\ [0.5em]
  & $E=1/(2\sqrt{2}) \simeq 0.353$, $Y_1=0$, $Y_0=1/\sqrt{2}$ \\ [0.5em]
   & $E=(\sqrt{2}+\sqrt{6})/4 \simeq 0.965$,  $Y_1=0$, $Y_0=1/\sqrt{2}$ \\ [0.5em]
     & $E=3\sqrt{2} /5 \simeq 0.848$, $Y_1 = 1/\sqrt{2}$, $Y_0 = 7/(5\sqrt{2}) \simeq 0.989$ \\ [0.5em]
     & $E=\sqrt{2} /5 \simeq 0.282$, $Y_1 = -1/\sqrt{2}$, $Y_0 =E/2 \simeq 0.141$ \\ [0.5em]
(5)  &   \\ [0.5em]
(6)  &   \\ [0.5em]
    \end{tabular}
  \end{center}
    \caption{Solutions to $(i)$, $1 \leq i \leq 6$, in \cite[subsection 4.2]{GM3} for $\kappa=4$.}
        \label{sol_kappa_2346}
\end{table}
\normalsize

\section{Polynomial Interpolation for the case $\kappa = 2$ in dimension 2, $5 \leq n \leq 9$}
\label{sectionk2}

In \cite[section 13]{GM3} we gave numerical evidence of bands of a.c.\ spectrum for bands $n=1,2,3,4$, for $\kappa=2$ in dimension 2. We briefly mention that :

$\bullet$ For the $5th$ band, we checked that $\Sigma = [2, 4, 6, 8, 10, 12, 14, 16, \alpha,\beta]$ is also valid for $(\alpha, \beta) = (18, 20), (20, 24), (22,24), (24,28), (24,32), (28,32)$ whereas it is not valid for $(\alpha, \beta) = (18,24)$, $(18,26), (18,28)$. 

$\bullet$ For the $6th$ band we checked that $\Sigma = [2, 4, 6, 8, 10, 12, 14, 16, 18, 20, \alpha,\beta]$ is valid for $(\alpha, \beta) = (22, 28), (22,30)$ but not valid for $(\alpha, \beta) = (22,24), (22,26)$. 

$\bullet$ For the $7th$ band we checked that $\Sigma = [2, 4, 6, 8, 10, 12, 14, 16, 18, 20, 22, 24, \alpha,\beta]$ is valid for $(\alpha, \beta) = (26, 28)$ but not valid for $(\alpha, \beta) = (28, 36)$. 

$\bullet$ For the $8th$ band we checked that $\Sigma = [2, 4, 6, 8, 10, 12, 14, 16, 18, 20, 22, 24, 26, 28, \alpha,\beta]$ is valid for $(\alpha, \beta) = (30,38), (30,40)$ but not valid for $(\alpha, \beta) = (30,32), (30,36)$. 

$\bullet$ For the $9th$ band we checked $\Sigma = [2, 4, 6, 8, 10, 12, 14, 16, 18, 20, 22, 24, 26, 28, 30, 32, 34, 36]$ is valid.

\section{Polynomial Interpolation for the case $\kappa = 3$ in dimension 2, $n=5$}
\label{sectionk3}

In \cite[section 14]{GM3} we gave numerical evidence of bands of a.c.\ spectrum for bands $n=1,2,3,4$, for $\kappa=3$ in dimension 2. We briefly give the numerical details for $n=5$, i.e.\ the $5^{th}$ band.  $X_3 = \E_5 /2$ and

\begin{equation*}
\label{sys1000001}
\begin{cases}
T_3(\E_5-1) = T_3(\E_5-X_1) & \\
T_3(X_1) = T_3(\E_5-X_2) & \\
T_3(X_2) = T_3(\E_5 /2) & \\
\end{cases}
\Rightarrow
\begin{cases}
4\big [ (\E_5-1)^2 + (\E_5-1)(\E_5-X_1) + (\E_5-X_1)^2 \big]  =3 & \\
4\big [ X_1^2 + X_1(\E_5-X_2)  +(\E_5-X_2) ^2 \big]  =3 & \\
4\big [ X_2^2 + X_2(\E_5 /2)  +(\E_5 /2) ^2 \big]  =3. & 
\end{cases}
\end{equation*}

Using Python's \textit{fsolve} we get $[\E_5, X_1, X_2] \simeq [1.112, \ 0.307, \ 0.441]$. For an analytical solution or perhaps more precise, expand the equations and get rid of the $X_1^2$ and $X_2^2$ terms. One ends up with $4 \left( 9 \E_5^2/4 - 4\E_5 X_1 +X_1 X_2 -3 \E_5 +1 +5 \E_5 X_2 /2  + X_1 \right) = 3.$
Then use the fact that 
$$X_1 = -\frac{1}{2} +\frac{3 \E_5}{2} - \frac{\sqrt{3}}{2} \sqrt{ \E_5 (2-\E_5)}, \quad X_2 = - \frac{\E_5}{4} +\frac{\sqrt{3}}{4} \sqrt{4-\E_5^2}.$$ (signs are chosen based on numerical solution).
We get an equation in $\E_5$ only, which is 
\begin{equation*}
\begin{aligned}
0 & = -(35 \E_5^2)/8  + \left(3 \E_5 - \sqrt{3} \sqrt{(2 - \E_5) \E_5} - 1 \right) \left((1/8) \sqrt{3} \sqrt{4 - \E_5^2} - \E_5 /8 \right) \\
& \quad + (5/8) \sqrt{3} \E_5 \sqrt{4 - \E_5^2} + 2 \sqrt{3} \E_5 \sqrt{(2 - \E_5) \E_5}  + \E_5 /2 - 1/2\sqrt{3} \sqrt{(2 - \E_5) \E_5} - 1/4.
\end{aligned}
\end{equation*}
$\E_5$ is a root of 
$$ mp(E) = 372775 E^8 - 750010 E^7 + 536359 E^6 - 270784 E^5 + 128593 E^4 - 36442 E^3 + 4333 E^2 - 220 E + 4.$$
and we suspect it's its minimal polynomial. $\E_5 \simeq 1.11207$. We are not aware of a closed formula for $\E_5$. It follows that $X_1 \simeq 0.30753$, $X_2 \simeq 0.44178$, $X_3 \simeq 0.55603$, $X_4 \simeq 0.67028$, $X_5 \simeq 0.80453$.


Next we suppose $\Sigma = [3,6,9,12,15,18,21,24,27,30]$. We fill the matrix $M$ with floats. Python says the solution to $M \rho =0$ is $[\rho_3, \rho_6, \rho_9,\rho_{12}, \rho_{15}, \rho_{18}, \rho_{21}, \rho_{24}, \rho_{27}, \rho_{30}]^T$ equals 
$$\simeq [1, 1.599931, 1.645307, 1.27734, 0.77838, 0.37292, 0.13741, 0.03703, 0.00657, 0.00058]^{T}.$$
Graphically it seems $x \mapsto G_{\kappa} ^E (x) >0$ for $E \in (\E_5, \E_4)$, $x \in [E-1,1]$. So $\Sigma$ is valid.

\section{Polynomial Interpolation for the case $\kappa = 4$ in dimension 2, $1 \leq n \leq 6$}
\label{sectionk4}

This section is the analogue of \cite[section 13]{GM3} and \cite[section 14]{GM3}, but for $\kappa=4$.

First we note that  
\begin{equation}
\label{T4}
T_4(x) = 8 x^4 - 8 x^2 + 1 \quad \text{and} \quad T_4(x) = T_4(y) \Leftrightarrow 8(x-y)(x+y) [x^2+y^2-1]=0.
\end{equation} 

Let $\E_n$ be the energy solutions of \cite[Proposition 7.1]{GM3} and \cite[Proposition 7.2]{GM3}, for $\kappa=4$. We have : $\E_n \in J_2(\kappa=4) = (2\cos(\pi/4), 1+ \cos(\pi/4)) = (\sqrt{2}, 1+\sqrt{2}/2)$, and
$$\E_n-1 = X_0 < X_1 < X_2 < ... < X_n <X_{n+1} := 1.$$

\subsection{$n=1$ : $1st$ band.} 
On the left, $\E_1 = 8/5$, $X_1 = 4/5$. On the right, $\E_0 = 1+\cos(\frac{\pi}{4}) = 1+\frac{\sqrt{2}}{2}$. We choose $\Sigma = [4,8]$ and get $\rho = [\rho_{4}, \rho_8 ]^{T} = [1, 625/1054]^{T}$. The function $x \mapsto G_{\kappa} ^E (x)$ is strictly positive for $E \in (\E_1, \E_0)$, $x \in [E-1,1]$. 

\begin{figure}[H]
  \centering
 \includegraphics[scale=0.178]{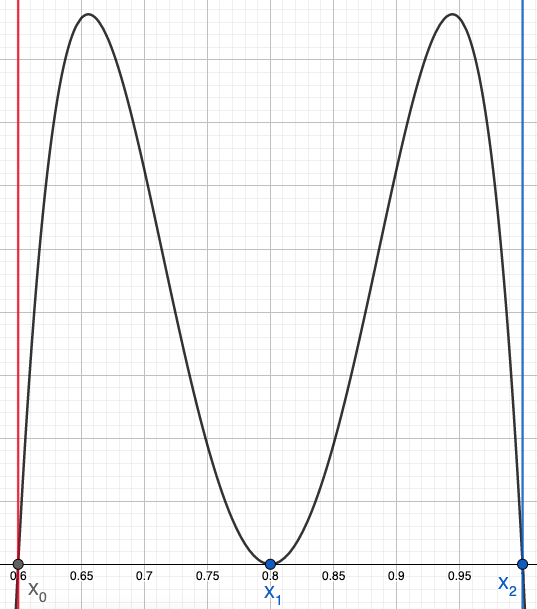}
  \includegraphics[scale=0.19]{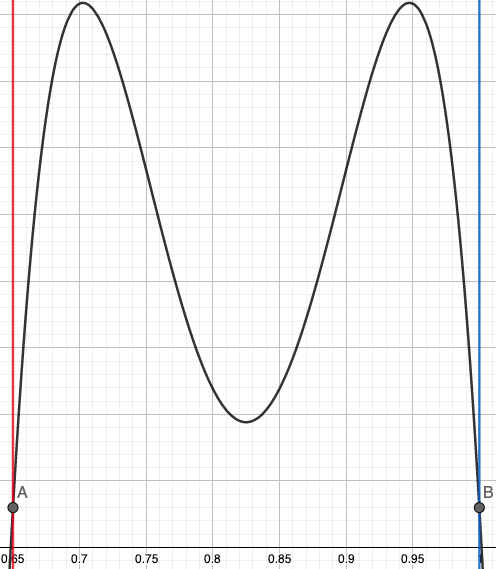}
  \includegraphics[scale=0.195]{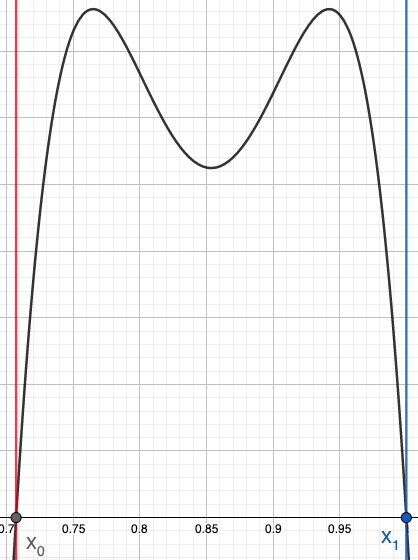}
\caption{$G_{\kappa=4} ^E (x)$. $x \in [E-1,1]$. From left to right: $E=\E_1$, $E=1.65$, $E=\E_0$. $G_{\kappa=4}^E (x) > 0$ in the middle picture, but not in the other two.}
\label{fig:test_k2_33oooi}
\end{figure}

\subsection{$n=2$ : $2nd$ band.} 
\begin{figure}[H]
  \centering
 \includegraphics[scale=0.17]{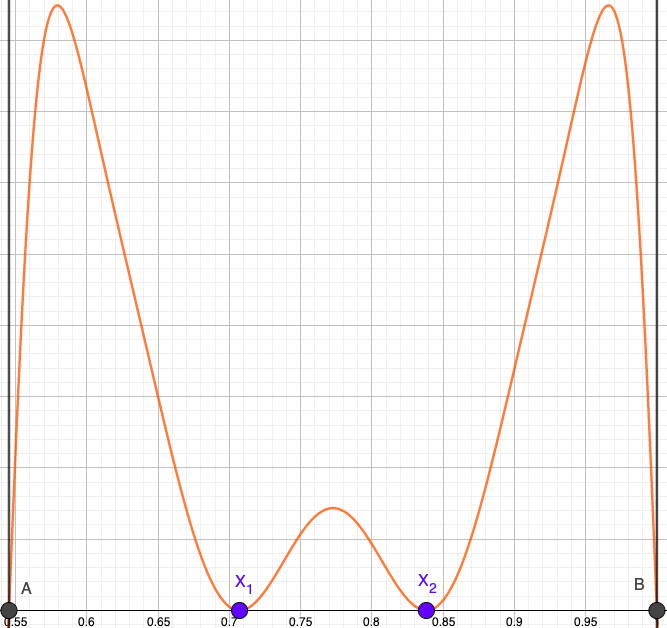}
  \includegraphics[scale=0.177]{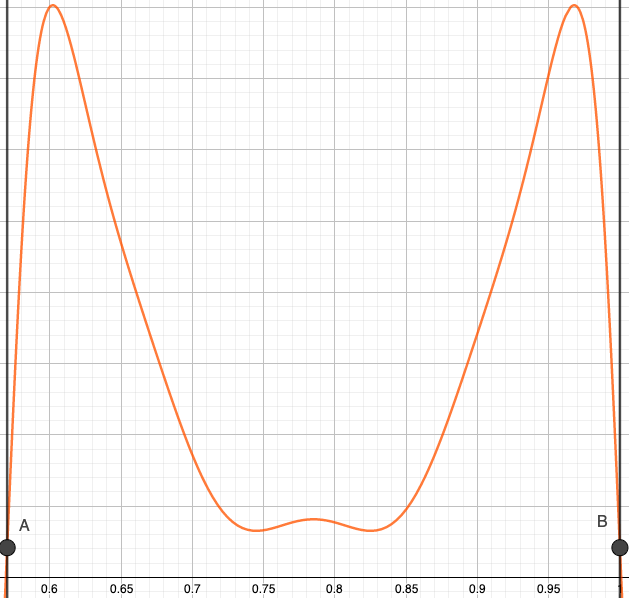}
  \includegraphics[scale=0.195]{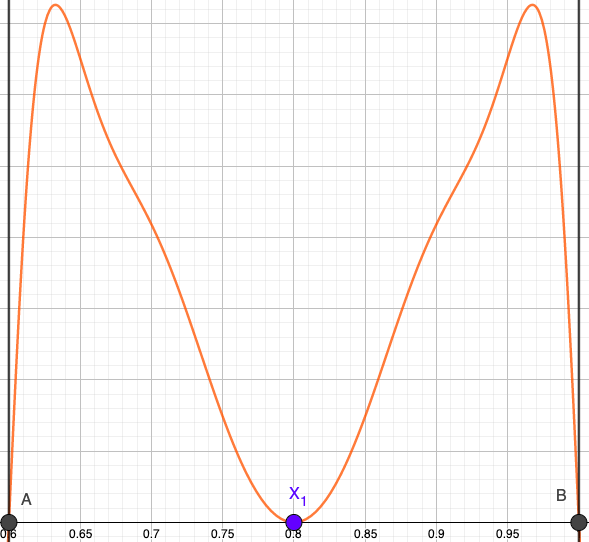}
\caption{$G_{\kappa=4} ^E (x)$. $x \in [E-1,1]$. From left to right: $E=\E_2$, $E=1.57$, $E=\E_1$. $G_{\kappa=4}^E (x) > 0$ in the middle picture, but not in the other two.}
\label{fig:test_k2_33333sish}
\end{figure}

$X_1 = \sqrt{2}/2$. Then we solve $T_4(\E_2-1) = T_4(\E_2 -\sqrt{2}/2)$ which has 2 solutions $\E_2 = (2 + \sqrt{2} \pm \sqrt{2 + 4 \sqrt{2}})/4$. Using the context, it must be that 
\begin{equation}
\label{E_L_k4_2}
\E_2 = \frac{2 + \sqrt{2} + \sqrt{2 + 4 \sqrt{2}}}{4} \simeq 1.54532.
\end{equation}
It follows that $X_2 = \E_2 - X_1 = (2 - \sqrt{2} + \sqrt{2 + 4 \sqrt{2}})/4$. We suppose $\Sigma = [4,8,12,24]$. The solution to $M \rho =0$ is $[\rho_4, \rho_8, \rho_{12}, \rho_{24}]^T \simeq [1, 0.81070, 0.21647, -0.06593]^{T}$. The function $x \mapsto G_{\kappa} ^E (x)$ is plotted in Figure \ref{fig:test_k2_33333sish} for some values $E \in [\E_2, \E_1]$.
Moreover it would appear that the combinations $\Sigma = [4,8,12,28], [4,8,12,32], [4,8,16,24]$ are equally valid, whereas $\Sigma = [4, 8, 12, 16], [4,8,12,20], [4,8,16,20]$ are not valid. The linear dependency between rows of $M$ seems to rely (at least partly) on the relations
$$U_{4j-1}(1/\sqrt{2}) =0, \quad \text{and} \quad U_{4j-1}(\sqrt{2}/\sqrt{3}) = - \sqrt{2} U_{4j-1}(1/\sqrt{3}), \quad j \in \N^*.$$

\subsection{$n=3$ : $3rd$ band.} 
\begin{figure}[H]
  \centering
 \includegraphics[scale=0.195]{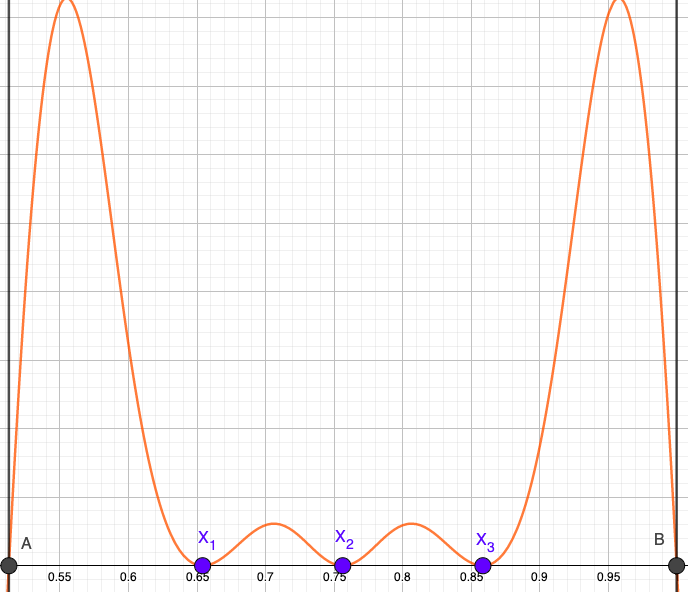}
  \includegraphics[scale=0.197]{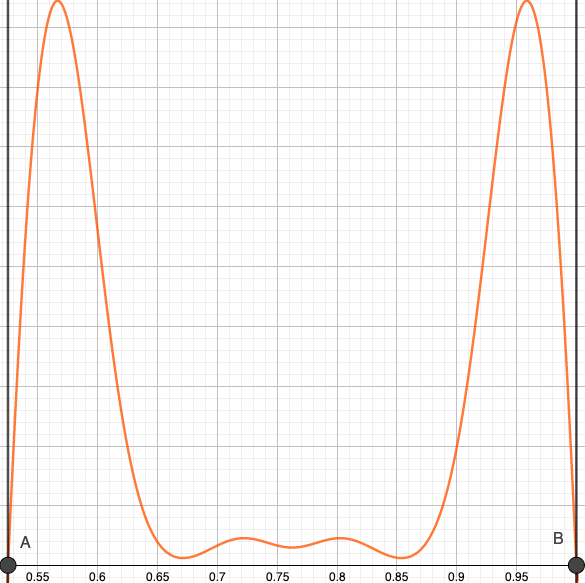}
  \includegraphics[scale=0.16]{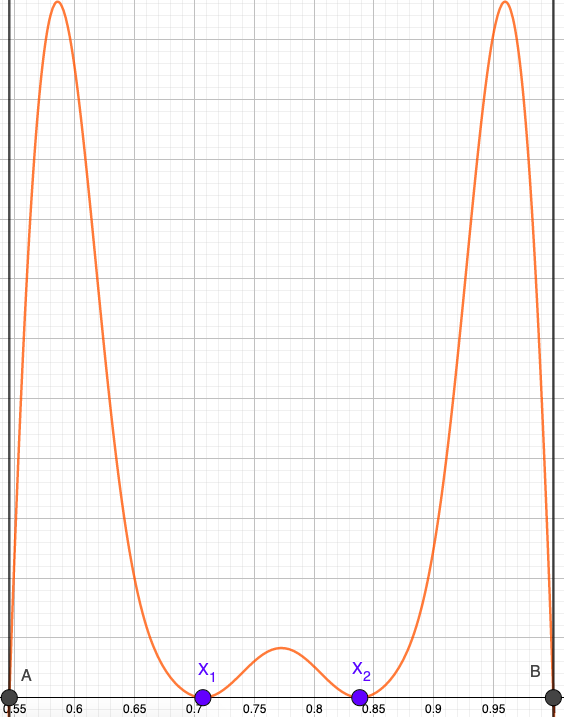}
\caption{$G_{\kappa=4} ^E (x)$. $x \in [E-1,1]$. From left to right: $E=\E_3$, $E=1.525$, $E=\E_2$. $G_{\kappa=4} ^E (x) > 0$ in the middle picture, but not in the other two. }
\label{fig:test_k2_33333sish2}
\end{figure}

$X_2 = \E_3/2$ and $T_4(\E_3-1) = T_4(\E_3-X_1)$, $T_4(X_1) = T_4(\E_3/2)$. Taking the context into account and applying \eqref{T4} leads to $(\E_3-1)^2 + (\E_3-X_1)^2 = 1$ and $X_1^2 +\E_3^2 /4 = 1$. So 
\begin{equation}
\label{E_L_k4_3}
\E_3 = \frac{28}{65} + \frac{2}{65 t}+ \frac{1}{2} \sqrt{ \frac{4256}{12675} + \frac{368}{195s} - \frac{16}{195} s + \frac{121088 t}{4225}} \simeq 1.51271, \quad t := \sqrt{\frac{3}{133 - \frac{1495}{s} + 65 s}}
\end{equation}
and $s :=  (629 + 48\sqrt{177})^{1/3}$. The minimal polynomial for $\E_3$ is $mp(E) = 65 E^4 - 112 E^3 + 56 E^2 - 64 E + 16$. Also $X_1 \simeq 0.65415$, $X_2 \simeq 0.75635$.
We choose the combination $\Sigma =[4, 8, 12, 16, 24, 36]$. The solution to $M \rho = 0$ is 
$$[\rho_4, \rho_8, \rho_{12}, \rho_{16}, \rho_{24}, \rho_{36}]^T \simeq [1, 1.18290, 0.68875, 0.18594, -0.00794, -0.00288]^{T}.$$
The function $G_{\kappa} ^E (x)$ is plotted in Figure \ref{fig:test_k2_33333sish2} for some values $E \in [\E_3, \E_2]$. The linear dependency between rows of $M$ seems to rely (at least partly) on the relation
$$ \left[U_{4j-1} (x), U_{4l-1} \left(\sqrt{1-x^2} \right) \right] = 0,\quad  \forall j,l \in \N^*,$$
where $[\cdot,\cdot ]$ is the Bezoutian defined by \cite[(2.3) of section 2]{GM3}.


The combination $\Sigma = [4,8,12,16,20,24]$ is equally valid and it gives 
$$[\rho_4, \rho_8, \rho_{12}, \rho_{16}, \rho_{20}, \rho_{24}]^T \simeq [1, 1.37002, 1.06973, 0.53992, 0.16964, 0.02655]^{T}.$$

\subsection{$n=4$ : $4th$ band.} 

$X_2 = \sqrt{2}/2$. Then we have to solve $(\E_4 - 1)^2 + (\E_4 - X_1)^2 =1$, $X_1^2 + (\E_4-\sqrt{2}/2)^2=1$. This leads to $1 - (\E_4 - \sqrt{2}/2)^2 = (\E_4 - \sqrt{\E_4(2-\E_4)})^2$. Thus $\E_4$ is the root of 
$$mp(E) = 5E^4 - (2\sqrt{2}+4) E^3 + (5 - 4 \sqrt{2}) E^2 + (\sqrt{2} -2) E + 1/4 = 0.$$
The closed form formula of $\E_4$ is horrendous; but $\E_4 \simeq 1.49137$. In turn, $X_1 \simeq 0.62042$. 

The simplest valid combination of indices we have found is $\Sigma = [4, 8, 12, 16, 20, 24, 28, 52]$. The solution to $M \rho = 0$ is
$\rho = [\rho_4, \rho_8, \rho_{12}, \rho_{16}, \rho_{20}, \rho_{24}, \rho_{28}, \rho_{52}]^{T} \simeq$
$$[1, 1.46864, 1.29941, 0.80098, 0.34657, 0.09808, 0.01417, 0.000030]^{T}.$$
Graphically $x \mapsto G_{\kappa} ^E (x)$ is strictly positive for $E \in (\E_4, \E_3)$, $x \in [E-1,1]$.

\subsection{$n=5$ : $5th$ band.} We have $X_3 = \E_5/2$, $(\E_5-1)^2 + (\E_5-X_1)^2 =1$, $X_1^2 + (\E_5-X_2)^2 = 1$, and $X_2^2 + (\E_5/2)^2 =1$. This leads to $1 - \left( \E_5 - \left(1 - \E_5^2/4 \right)^{1/2} \right)^2 = \left( \E_5 - \E_5^{1/2}(2-\E_5)^{1/2} \right)^2.$
The solution is 
\begin{equation}
\label{E_L_k4_5}
\E_5 = \frac{16}{1275} \left(42 + \frac{w}{3^{2/3}} - \frac{32533 w^{-1}}{3^{1/3}} \right) \simeq 1.47650, \quad w := (13025367 + 208250 \sqrt{6294})^{1/3}
\end{equation}
and its minimal polynomial is $mp(E) = -32768 + 22784 E - 6048 E^2 + 3825 E^3.$ In turn $X_1 \simeq 0.59734$, $X_2 \simeq 0.67452$, $X_3 \simeq 0.73825$, $X_4 \simeq 0.80198$ and $X_5 \simeq 0.87916$. We choose $\Sigma = [4, 8, 12, 16, 20, 24, 28, 32, 36, 40]$. $M \rho = 0$ gives
$[\rho_4, \rho_8, \rho_{12}, \rho_{16}, \rho_{20}, \rho_{24}, \rho_{28}, \rho_{32}, \rho_{36}, \rho_{40}]^{T} \simeq$
$$[1, 1.58691, 1.60962, 1.22543, 0.72779, 0.33759,  0.11956, 0.03071, 0.00514, 0.00042]^{T}.$$
It appears graphically that $x \mapsto G_{\kappa} ^E (x)$ is strictly positive for $E \in (\E_5, \E_4)$, $x \in [E-1,1]$.

\subsection{$n=6$ : $6th$ band.} $X_3 = \sqrt{2}/2$, $(\E_6-1)^2 + (\E_6-X_1)^2 =1$, $X_1^2 + (\E_6-X_2)^2 = 1$ and $X_2^2 + (\E_6-\sqrt{2}/2)^2 =1$.  Performing elementary operations shows that $\E_6 \simeq 1.46568$ is a root of the polynomial
\footnotesize
$$mp(E) = 384 E^6 + (-448 - 224 \sqrt{2}) E^5 + (200 + 128 \sqrt{2}) E^4 + (-72 - 52\sqrt{2}) E^3 + (16 + 12 \sqrt{2}) E^2 + (-2 - \sqrt{2}) E + 1/8.$$
\normalsize
In turn $X_1 \simeq 0.58072$, $X_2 \simeq 0.65158$, $X_3 \simeq 0.70710$, $X_4 \simeq 0.75857$, $X_5 \simeq 0.81409$, $X_6 \simeq 0.88495$.

\section{Application of Polynomial Interpolation for the first band ($n=1$), $2 \leq \kappa \leq 9$} 
\label{Appn=1}
This section is about the first band ($n=1$) in the interval $J_2(\kappa)$.
For the $1st$ band we suppose only 2 terms in the linear combination \cite[(1.6) of Introduction]{GM3} are required. In other words, we suppose $\Sigma = [j_1 \kappa, j_2 \kappa]$. According to \cite[section 12]{GM3} we look to solve the system
\begin{equation}
\label{system1p}
\begin{cases}
g _{j_1\kappa} ^{\E_1} (x) + \rho_{j_2\kappa} \cdot g^{\E_1} _{j_2\kappa} (x)  = 0, & \text{at} \ x = X_0 = \E_1-1, \\[0.3em]
g^{\E_1} _{j_1\kappa} (x) + \rho_{j_2\kappa} \cdot g^{\E_1} _{j_2\kappa} (x)  = 0, & \text{at} \ x = X_1, \\[0.3em]
\frac{d}{dx}g^{\E_1} _{j_1\kappa} (x) + \rho_{j_2\kappa} \cdot  \frac{d}{dx} g^{\E_1} _{j_2\kappa} (x)  = 0, & \text{at} \ x = X_1, \\[0.3em]
g^{\E_0} _{j_1\kappa} (x) + \rho_{j_2\kappa} \cdot  g^{\E_0} _{j_2\kappa} (x)  = 0, & \text{at} \ \E_0 = 1+\cos(\pi / \kappa), x = \E_0 -1. 
\end{cases}
\end{equation}
Note we are assuming $\rho_{j_1 \kappa} =1$. By \cite[Remark 5.2]{GM3}, $X_1 = \E_1/2$. So the $3rd$ line of \eqref{system1p} is trivially true, by \cite[Lemma 3.4]{GM3}. The $4th$ line is also always true, by \cite[Lemma 1.4]{GM3}. Assuming $\E_1 \neq 0, \pm2$, the system \eqref{system1p} is equivalent to 
\begin{equation}
\label{system2p}
\begin{cases} 
[U_{j_1\kappa-1}(\E_1-1) , U_{j_2\kappa-1}(\E_1/2) ]  =0, \\
X_1 = \E_1/2, \\
\rho_{j_2\kappa} = - U_{j_1\kappa-1}(\E_1-1) / U_{j_2\kappa-1}(\E_1-1) = - U_{j_1\kappa-1}(\E_1/2) / U_{j_2\kappa-1}(\E_1/2) .
\end{cases}
\end{equation}
Let $X_2 := 1$, $X_0 = \E_1-X_2 = \E_1 -1$. From the $3rd$ line of \eqref{system2p} we see that we want $U_{j \kappa-1}(X_0) \neq 0$ and $U_{j \kappa-1}(X_1) \neq 0$, in particular $U_{\kappa-1}(X_0), U_{\kappa-1}(X_1) \neq 0$. Thus, by \cite[Corollary 2.3]{GM3} we see that \eqref{system2p} is equivalent to 
\begin{equation}
\label{system44p}
\begin{cases} 
X_2 = 1, \\
X_0 = \E_1 - X_2, \\
X_1 = \E_1/2, \\
T_{\kappa}(X_0) = T_{\kappa}(X_1),\\
\rho_{j_2\kappa} = - U_{j_1\kappa-1}(\E_1-1) / U_{j_2\kappa-1}(\E_1-1) = - U_{j_1\kappa-1}(\E_1/2) / U_{j_2\kappa-1}(\E_1/2) .
\end{cases}
\end{equation}

This is in agreement with \cite[(5.1) of section 5]{GM3} for $n=1$. In practice, our assessment is that it is best to simply take $\Sigma = \{ j_1\kappa, j_2 \kappa \} = \{\kappa, 2\kappa \}$. Table \ref{table_look} lists the results for $\kappa \in \{2,...,9\}$.

\begin{table}[H]
  \begin{center}
    \begin{tabular}{c|c|c|c} 
      $\kappa$ & Solution to \eqref{system44p} & $\kappa$ & Solution to \eqref{system44p} \\ [0.5em]
      \hline
2 & $\E_1 = \frac{2}{3}$  & 3 & $\E_1 = \frac{5+3 \sqrt{2}}{7} \simeq 1.3203$ \\[0.5em]
2 & $\rho_{4} = \frac{9}{14} \simeq 0.6428$  & 3 & $\rho_{6} = \frac{170 - 81 \sqrt{2}}{92} \simeq  0.6027$ \\[0.5em]
4 &  $\E_1 =  \frac{8}{5}$ & 5 & $\E_1 = \frac{49 + 12 \sqrt{10} + \sqrt{5 (49 + 12 \sqrt{10})}}{62} \simeq 1.7386$ \\[0.5em]
4 & $\rho_{8} = \frac{625}{1054} \simeq 0.5929$ & 5 & $\rho_{10} = \frac{-4000073 + 2667375 \sqrt{2} + 225 \sqrt{5 (25786331 - 6299370 \sqrt{2})}}{3122396} \simeq 0.5889$ \\[0.5em]
6 & $\E_1 = 1 + \sqrt{\frac{2}{3}} \simeq 1.8164$ & 7 & $\E_1 \simeq 1.8642$ (closed form unknown to us)  \\[0.5em]
6 & $\rho_{12} =  \frac{27}{46} \simeq 0.5869$ & 7 & $\rho_{14} \simeq 0.5857$ \\[0.5em]
8 & $\E_1 = \frac{16 + 3 \sqrt{2} + 2 \sqrt{26 + 7 \sqrt{2}}}{17}$ & 9 &  $\E_1 \simeq 1.9173$ (closed form unknown to us) \\[0.5em]
8 & $\E_1 \simeq 1.8956$ &  &  \\[0.5em]
8 & $\rho_{16} \simeq 0.5850$ & 9 & $\rho_{18} \simeq 0.5844$
    \end{tabular}
  \end{center}
    \caption{Solutions to \eqref{system44p}. $\E_1 = 1st$ band's left endpoint.}
        \label{table_look}
\end{table}

\section{Energy thresholds in $(0,1/2)$ for $\kappa=3$ in dimension 2}
\label{section work in progress}

Fix $\kappa=3$, in dimension 2. We construct graphically a bunch of thresholds in $(0,1/2)$. Recall that by \cite[Lemma 5.3]{GM3}, $\mathfrak{T}_{n, \kappa} = - \mathfrak{T}_{n, \kappa}$, or by \cite[Lemma 3.2]{GM3}, $\boldsymbol{\Theta}_{\kappa} (\Delta) = - \boldsymbol{\Theta}_{\kappa} (\Delta)$. So the negative thresholds listed below have a positive counterpart.

\begin{figure}[H]
  \centering
 \includegraphics[scale=0.165]{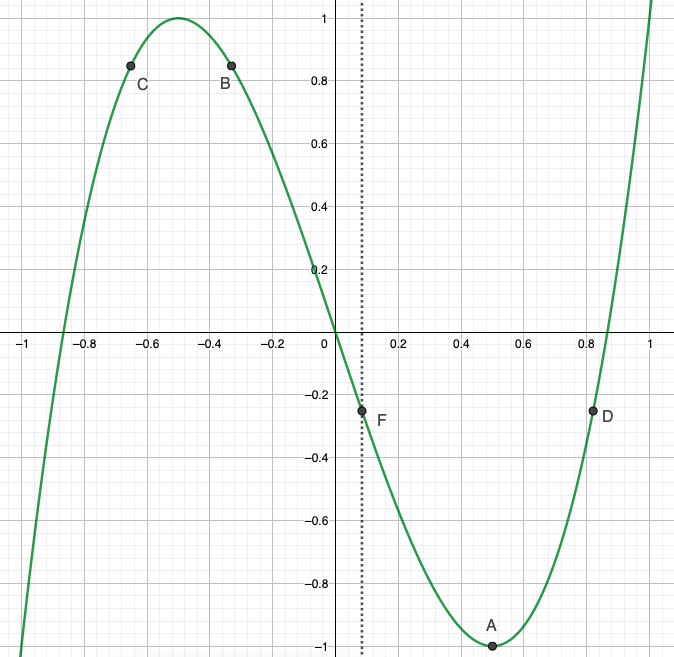}
  \includegraphics[scale=0.165]{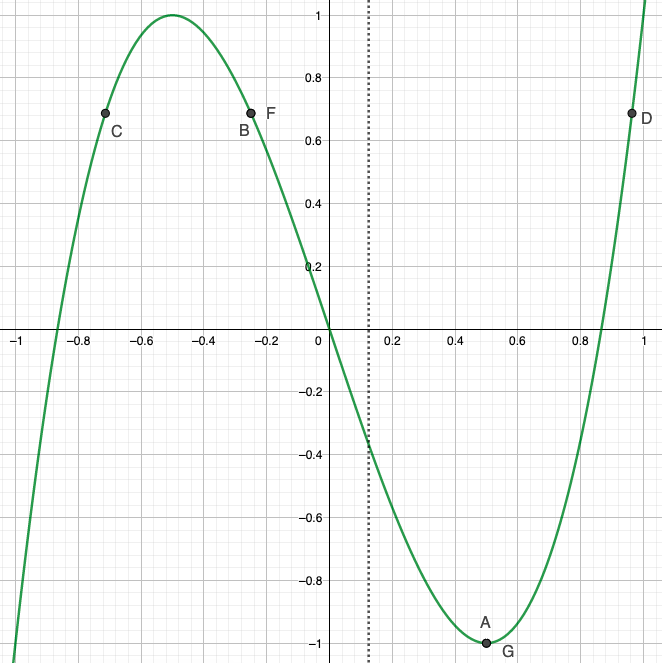}
  \includegraphics[scale=0.18]{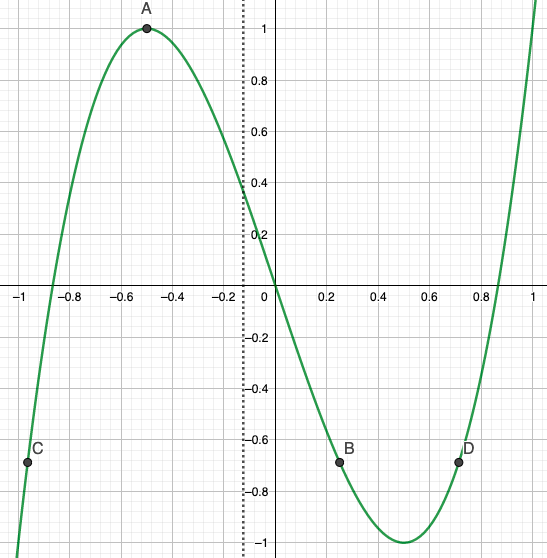}
  \includegraphics[scale=0.165]{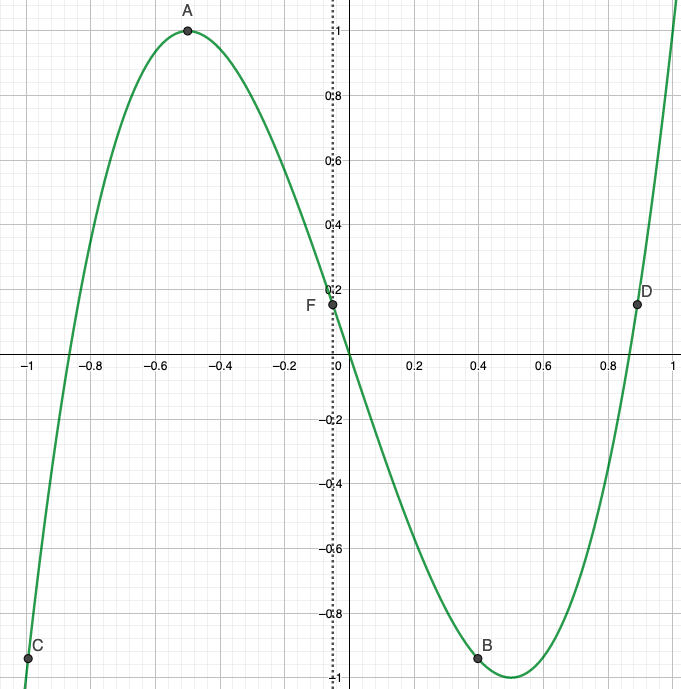}
   \includegraphics[scale=0.165]{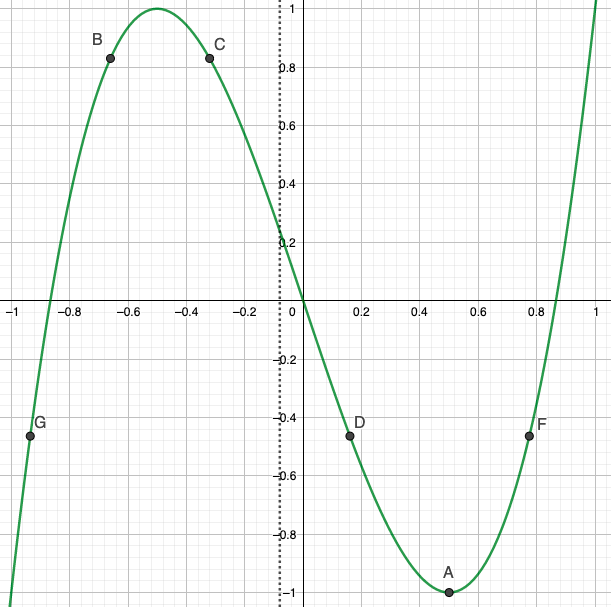}
      \includegraphics[scale=0.165]{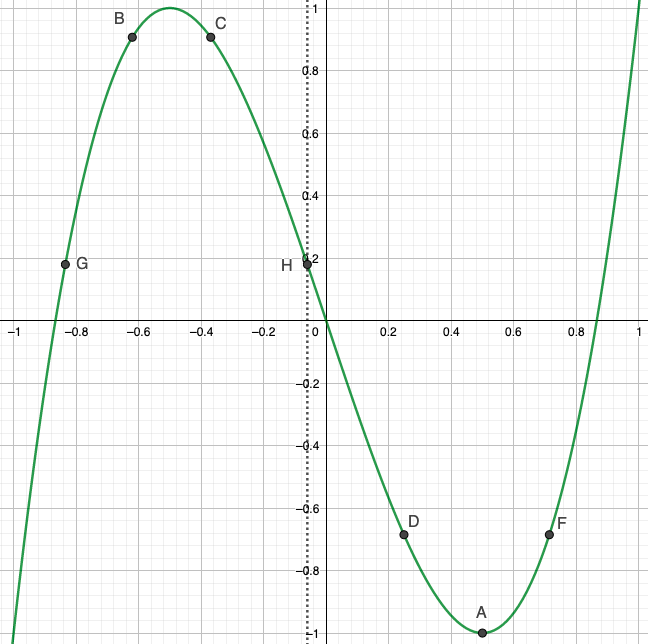}
       \includegraphics[scale=0.165]{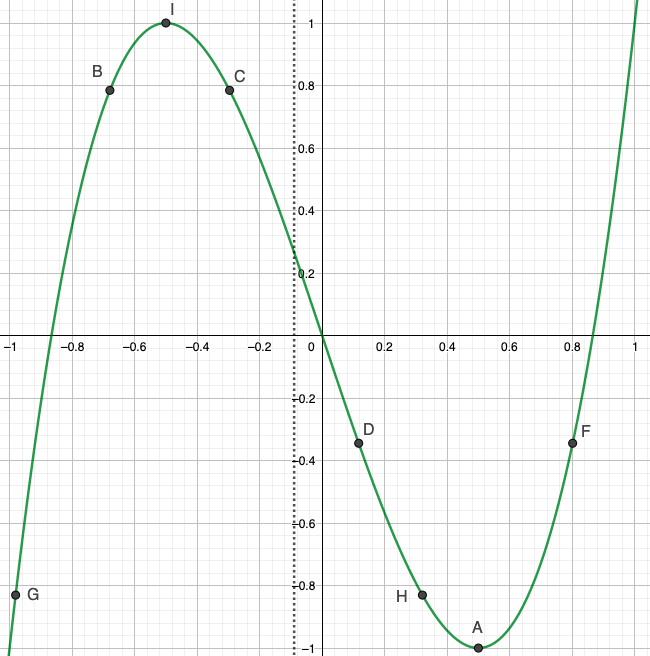}
       \includegraphics[scale=0.165]{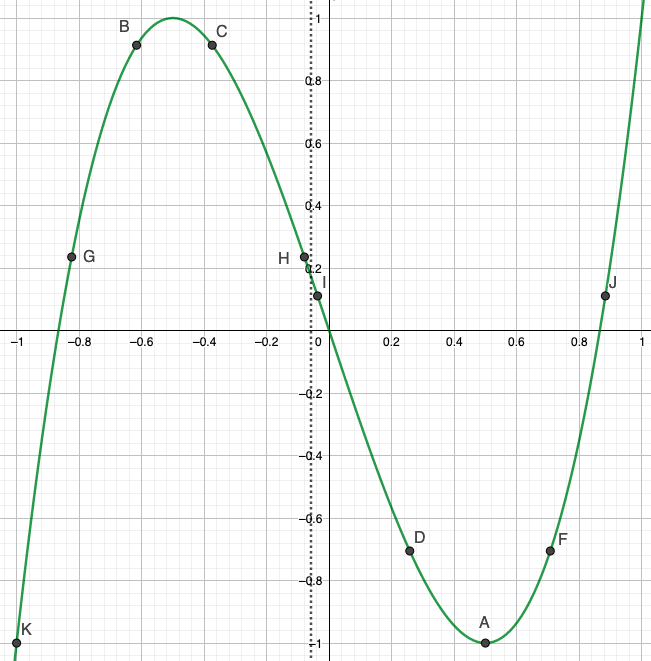}
         \includegraphics[scale=0.165]{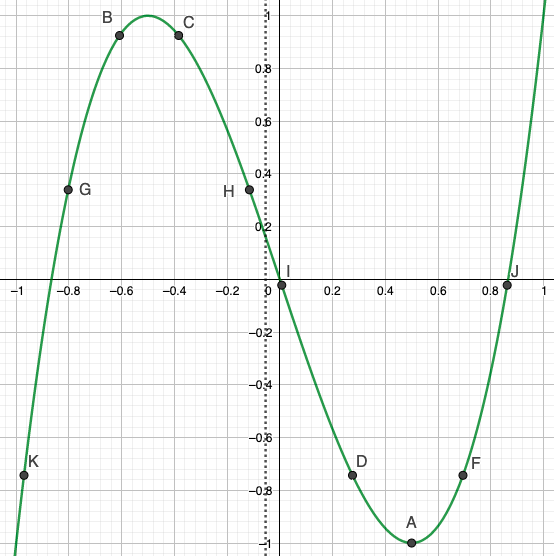}
           \includegraphics[scale=0.165]{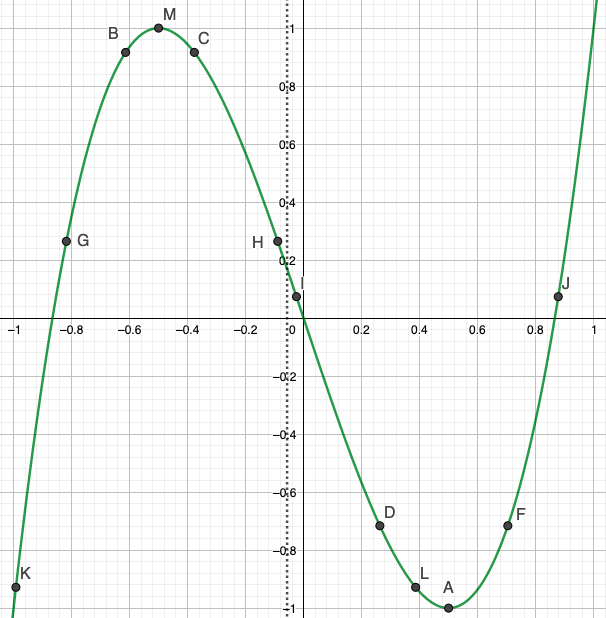}
           \includegraphics[scale=0.165]{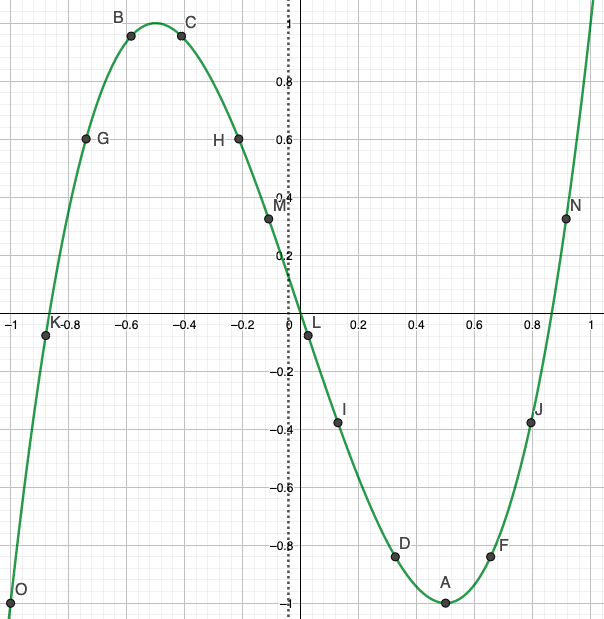}
            \includegraphics[scale=0.165]{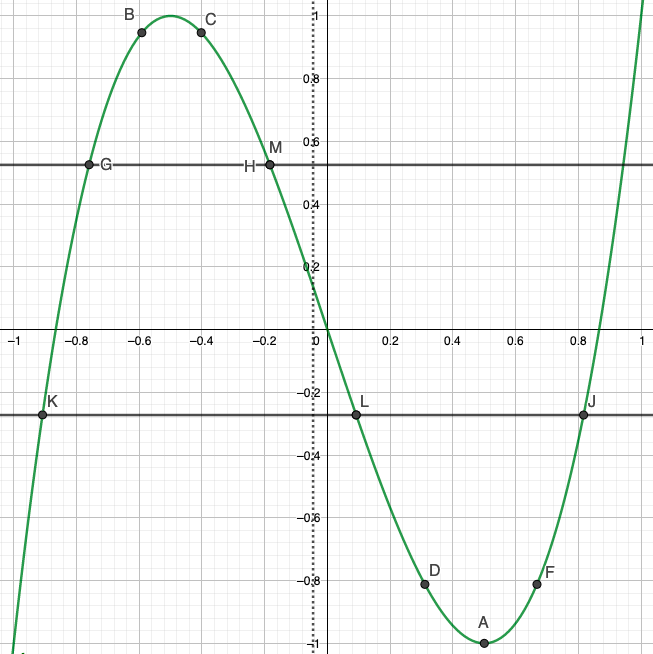}
           \includegraphics[scale=0.165]{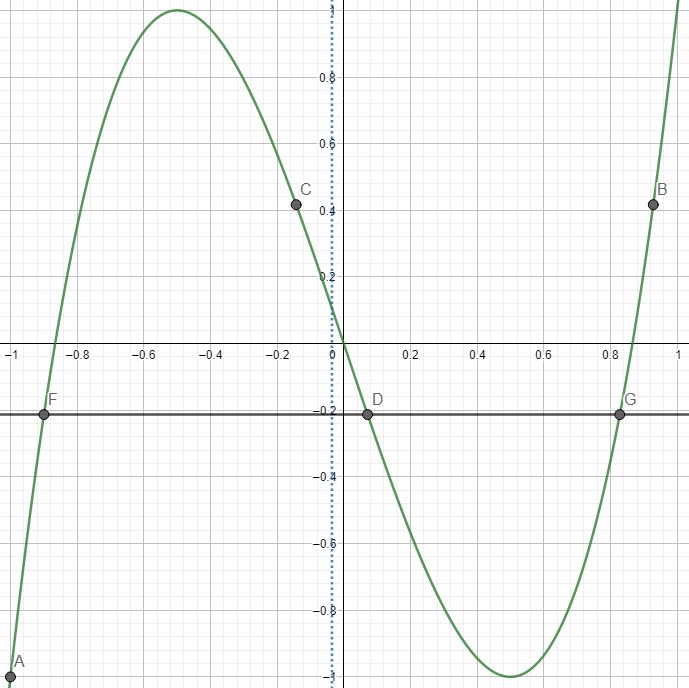}
           \includegraphics[scale=0.165]{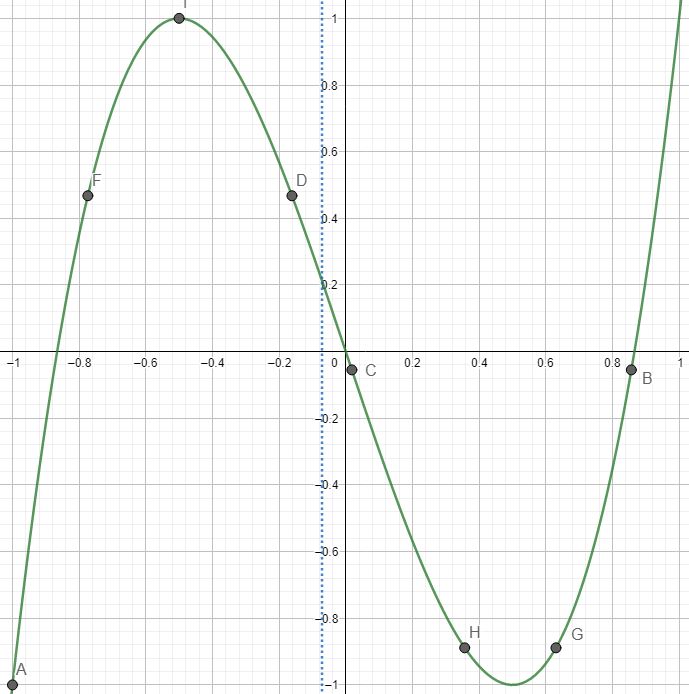}
            \includegraphics[scale=0.165]{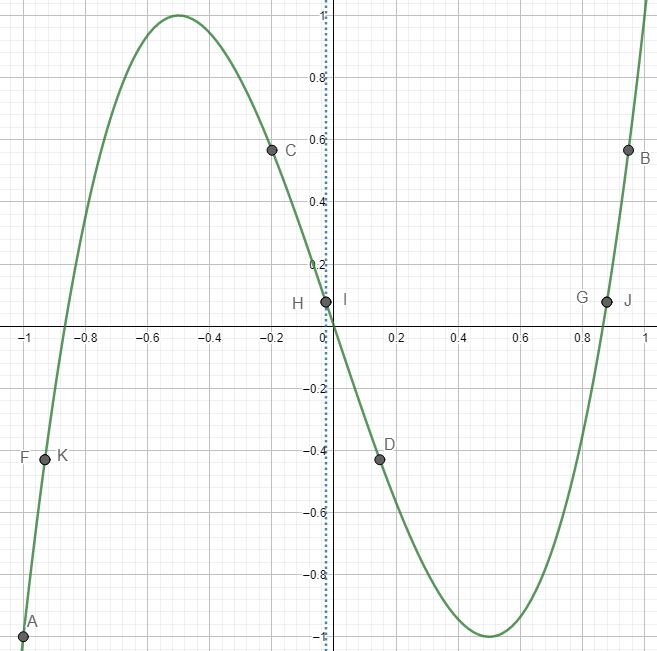}
            \includegraphics[scale=0.165]{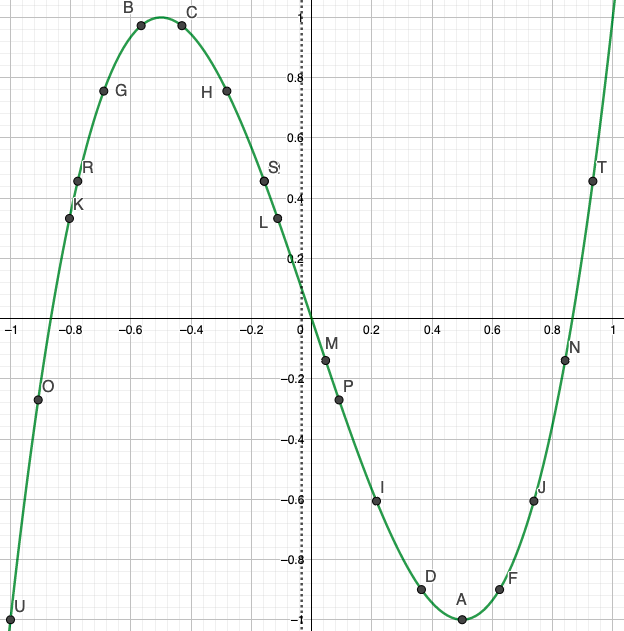}
            \includegraphics[scale=0.165]{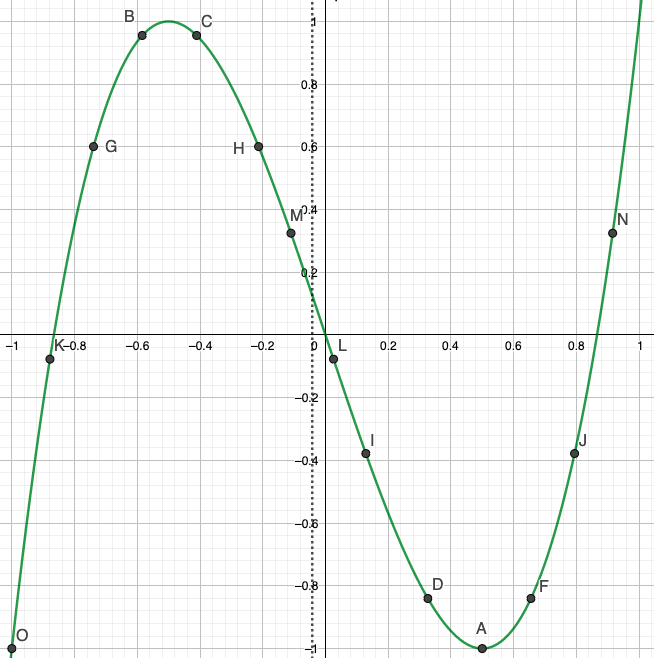}
            \includegraphics[scale=0.165]{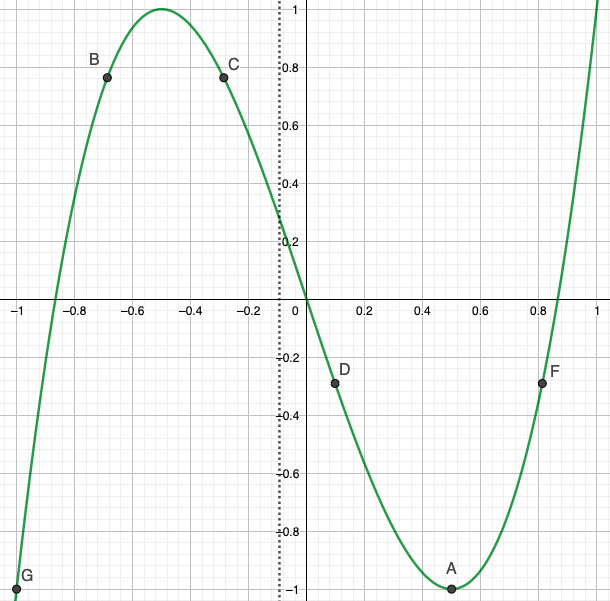}
            \includegraphics[scale=0.165]{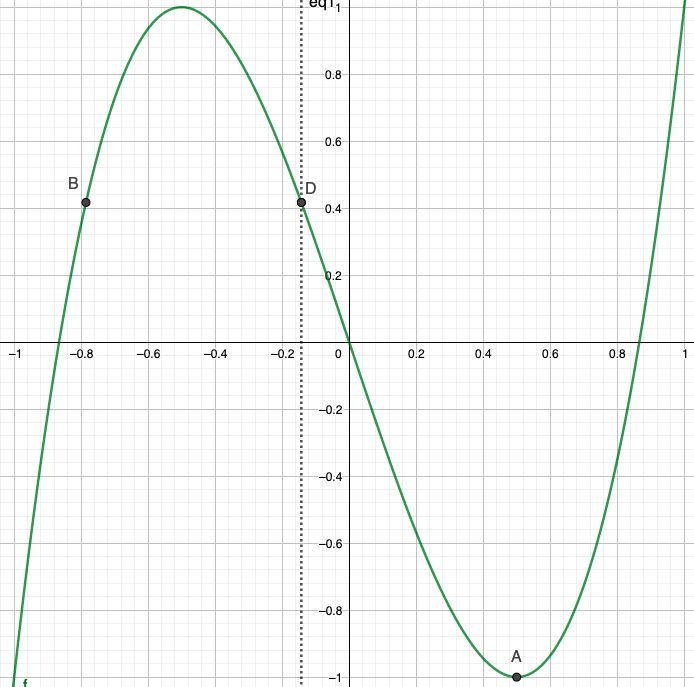}
            \includegraphics[scale=0.165]{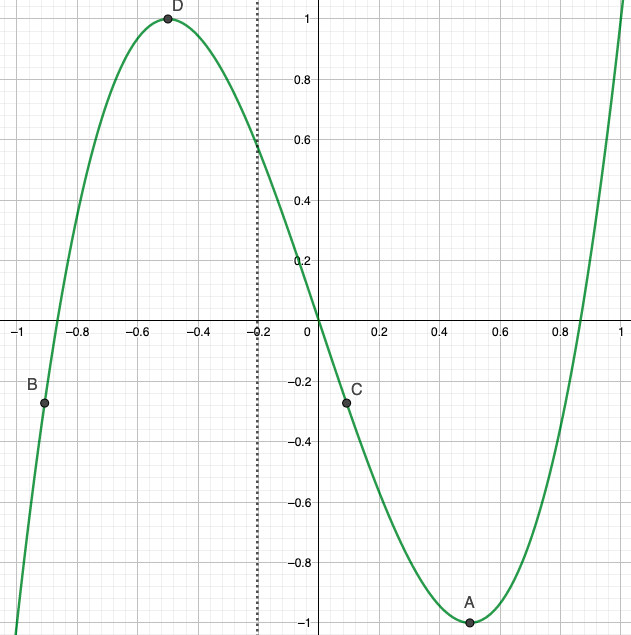}
             \includegraphics[scale=0.165]{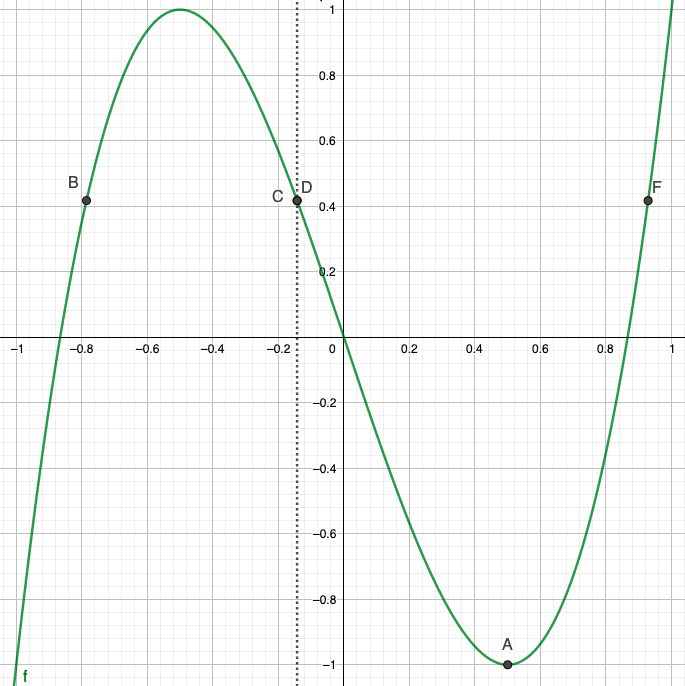}
              \includegraphics[scale=0.165]{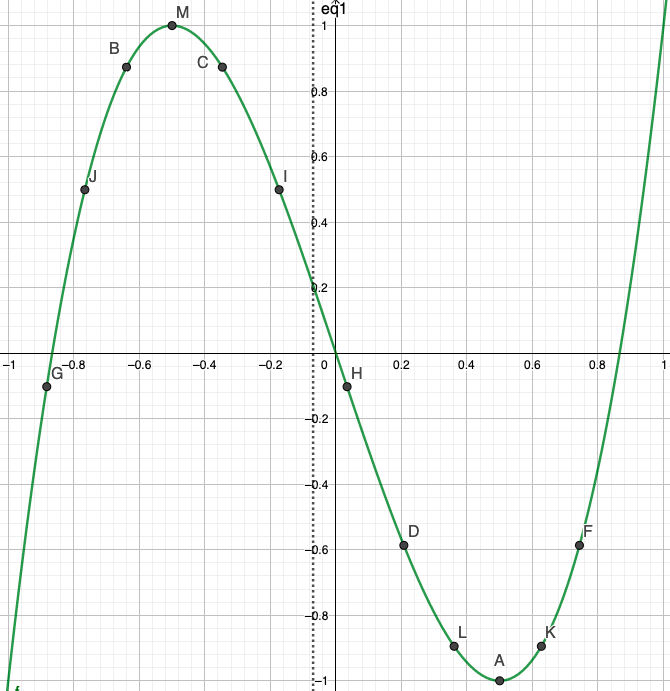}
              \includegraphics[scale=0.165]{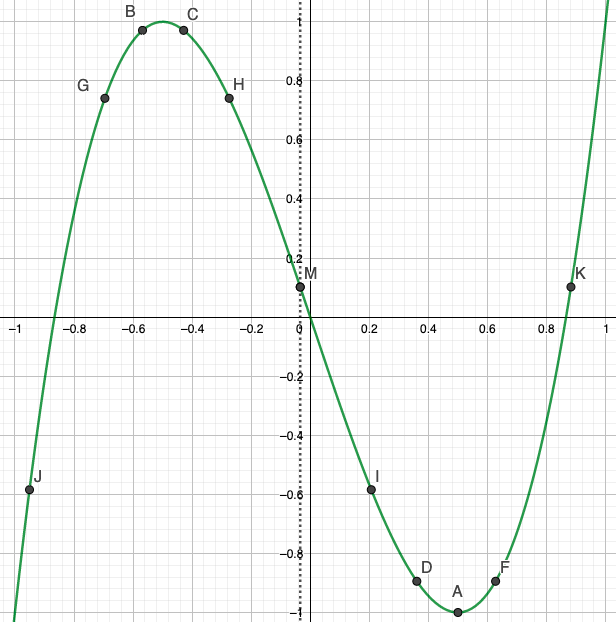}
              \includegraphics[scale=0.165]{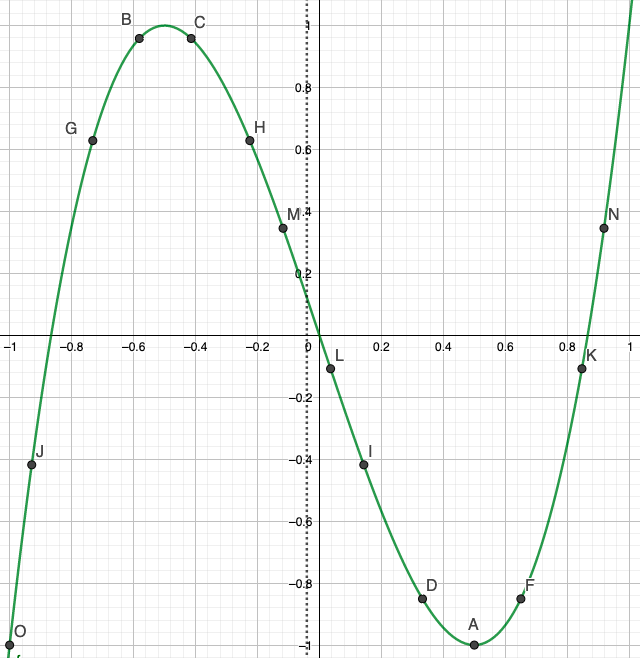}
\caption{}
\label{fig:work_in_progress}
\end{figure}

\begin{figure}[H]
  \centering
               \includegraphics[scale=0.165]{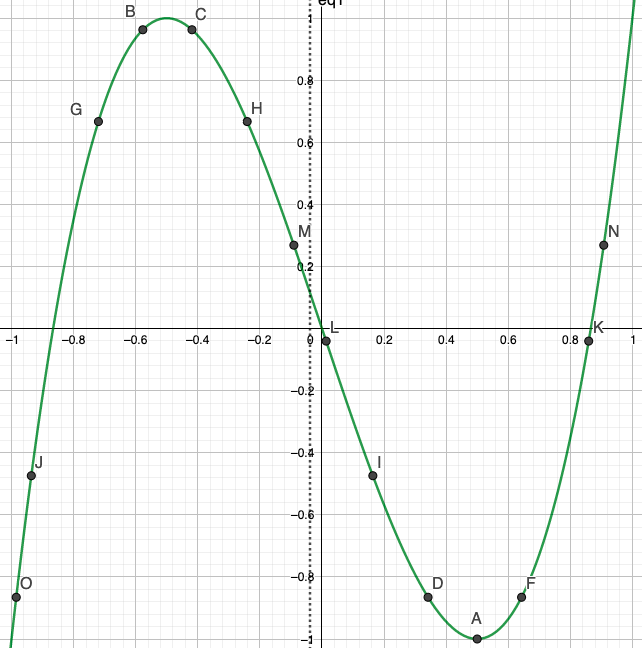}
               \includegraphics[scale=0.165]{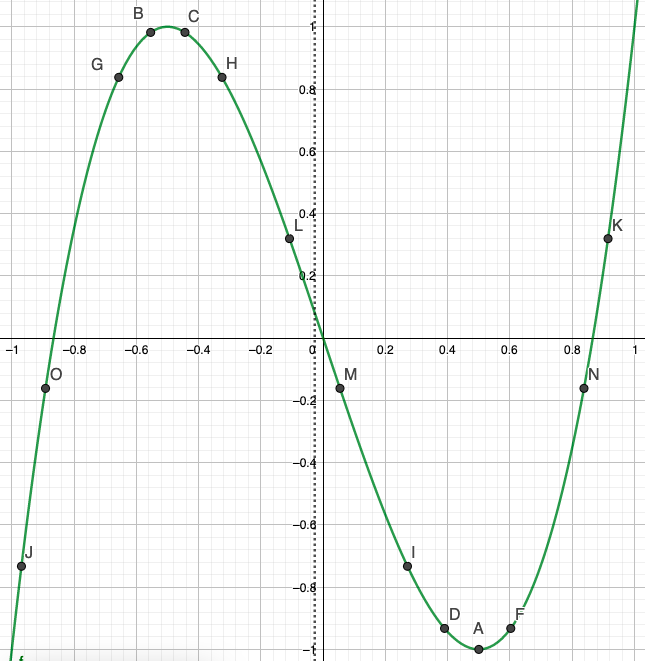}
                \includegraphics[scale=0.165]{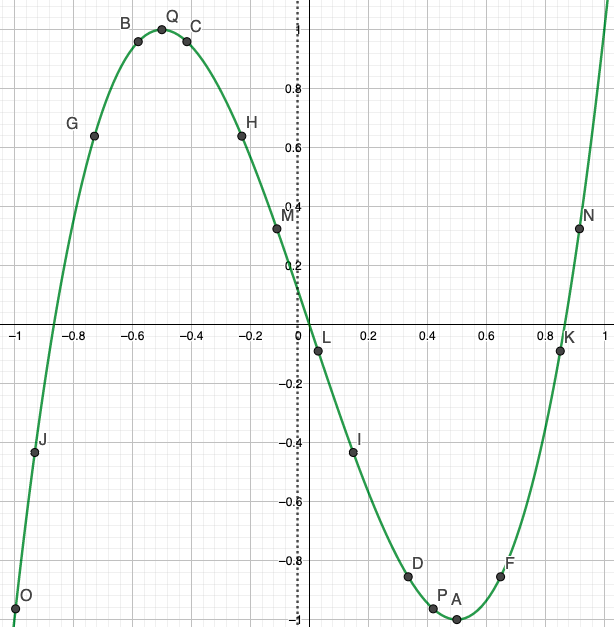}
               \includegraphics[scale=0.165]{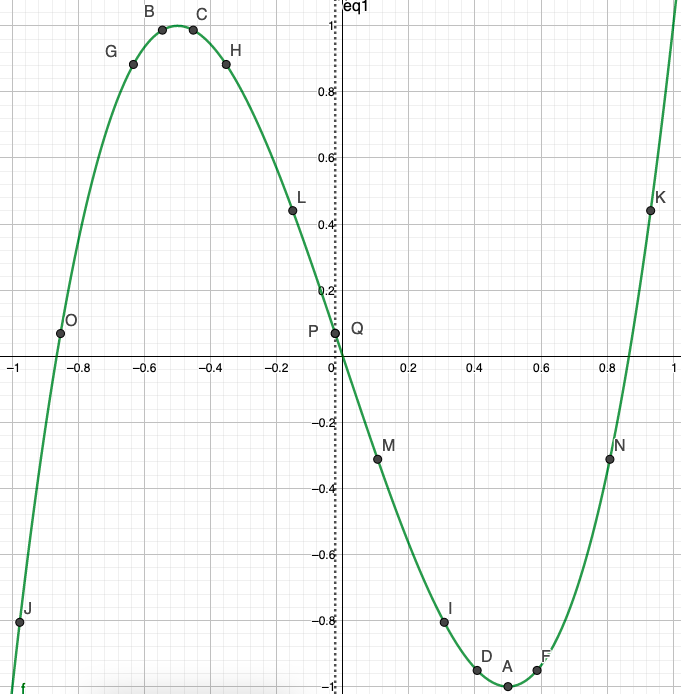}
               \includegraphics[scale=0.165]{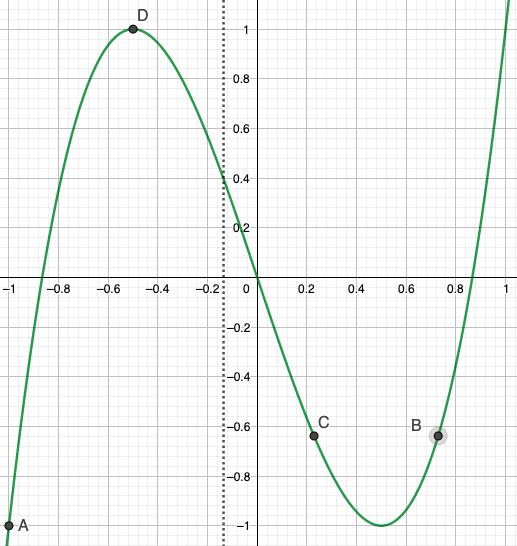}
                \includegraphics[scale=0.165]{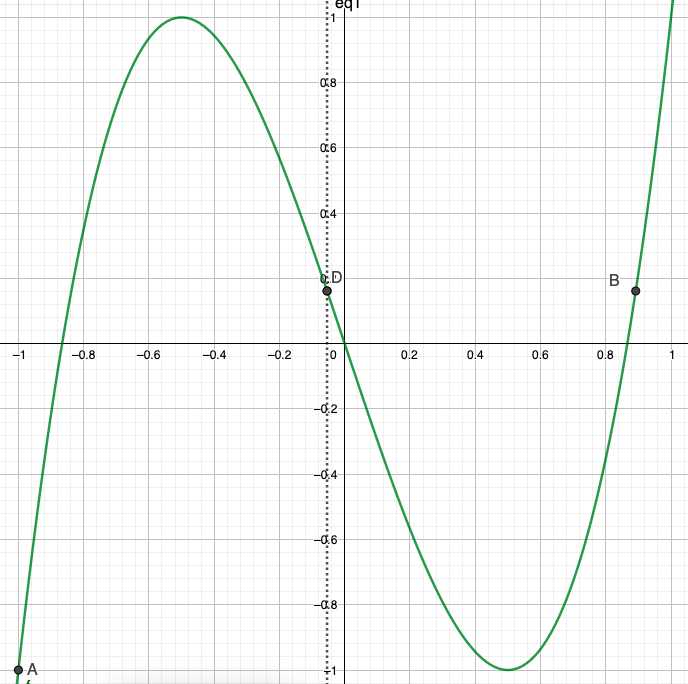}
                 \includegraphics[scale=0.165]{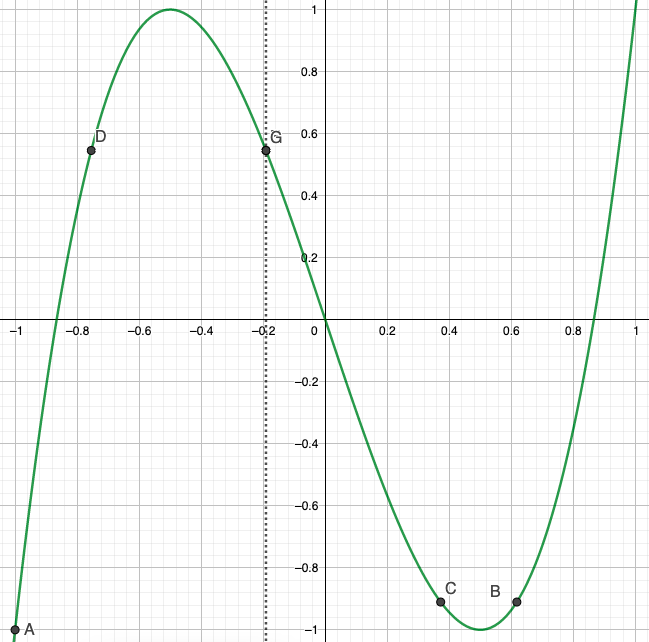}
                  \includegraphics[scale=0.165]{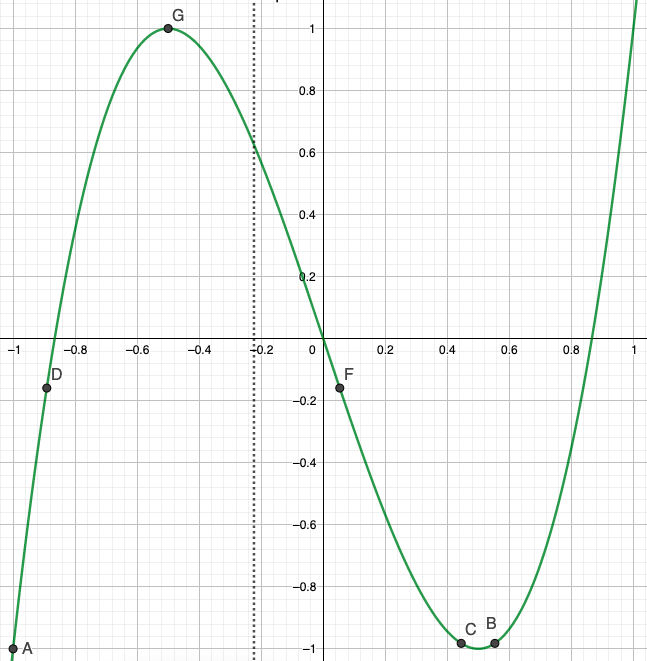}
                  \includegraphics[scale=0.165]{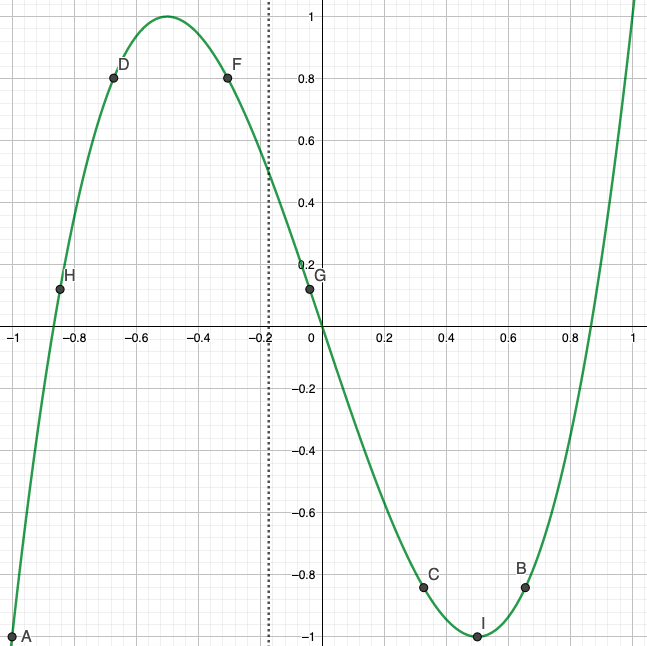}
                    \includegraphics[scale=0.165]{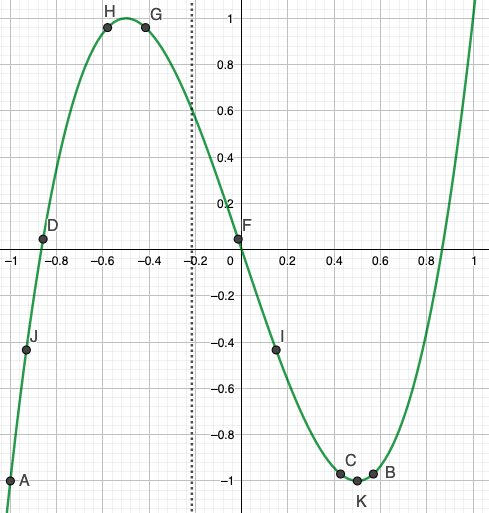}
                     \includegraphics[scale=0.165]{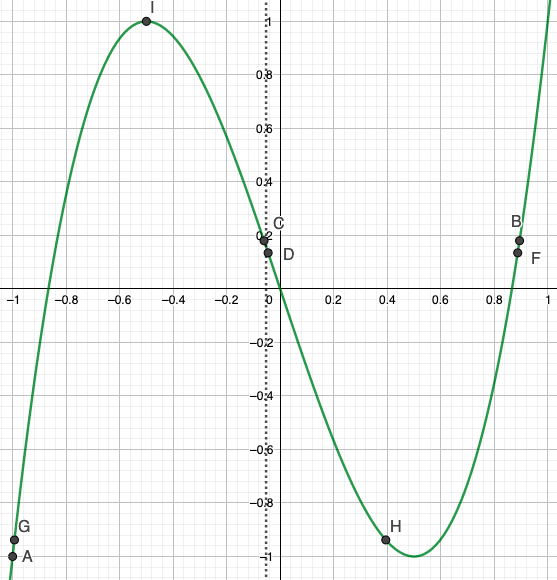}
                     \includegraphics[scale=0.165]{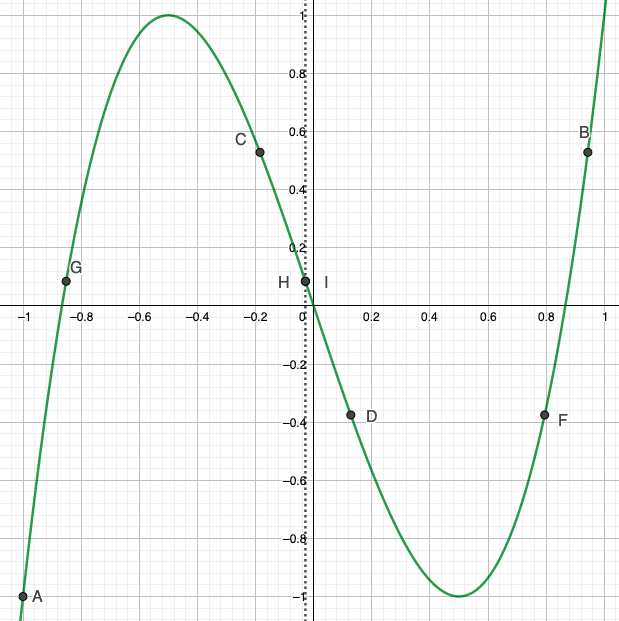}
                      \includegraphics[scale=0.165]{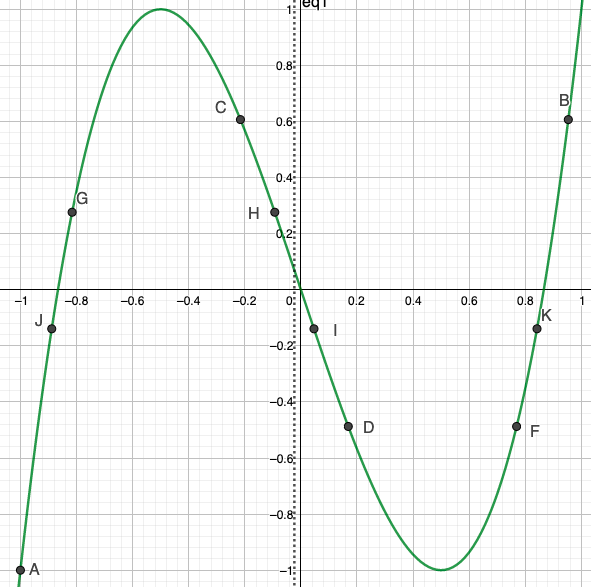}
                       \includegraphics[scale=0.165]{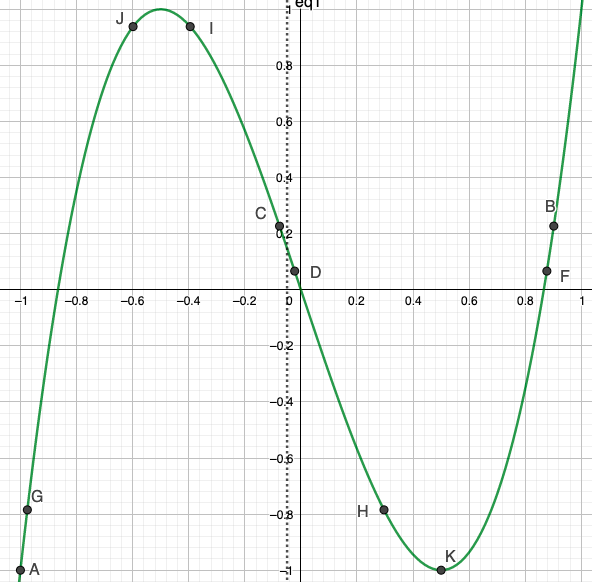}
                       \includegraphics[scale=0.165]{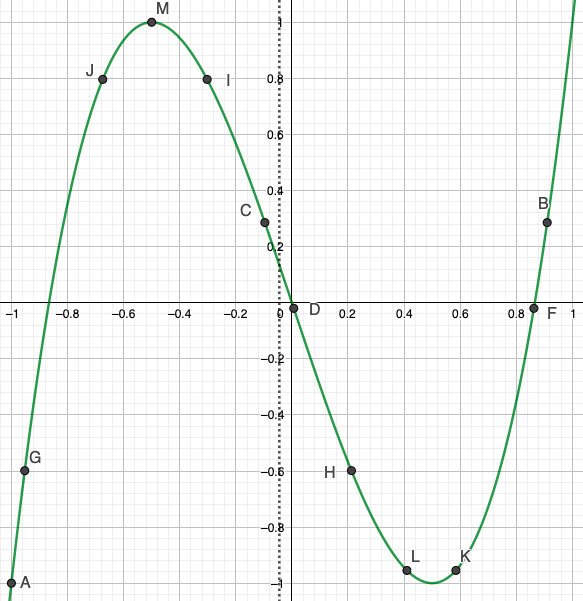}
                       \includegraphics[scale=0.15]{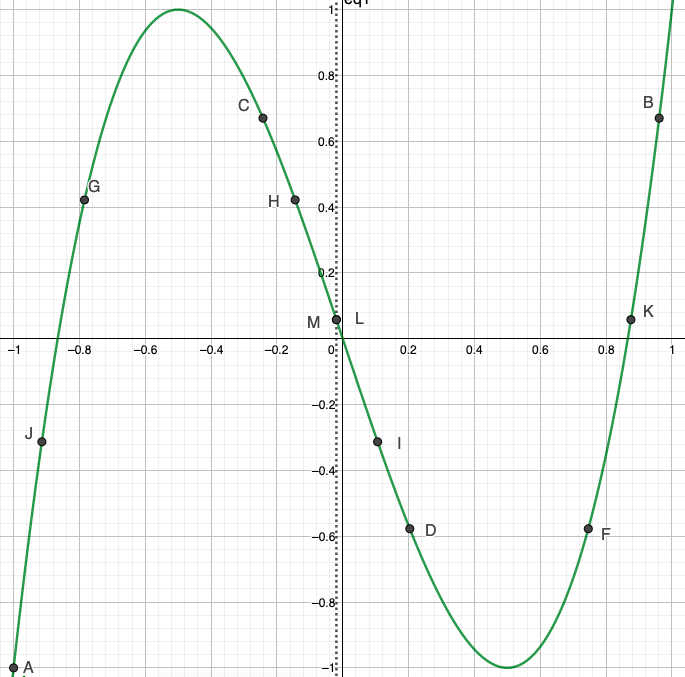}
                        \includegraphics[scale=0.15]{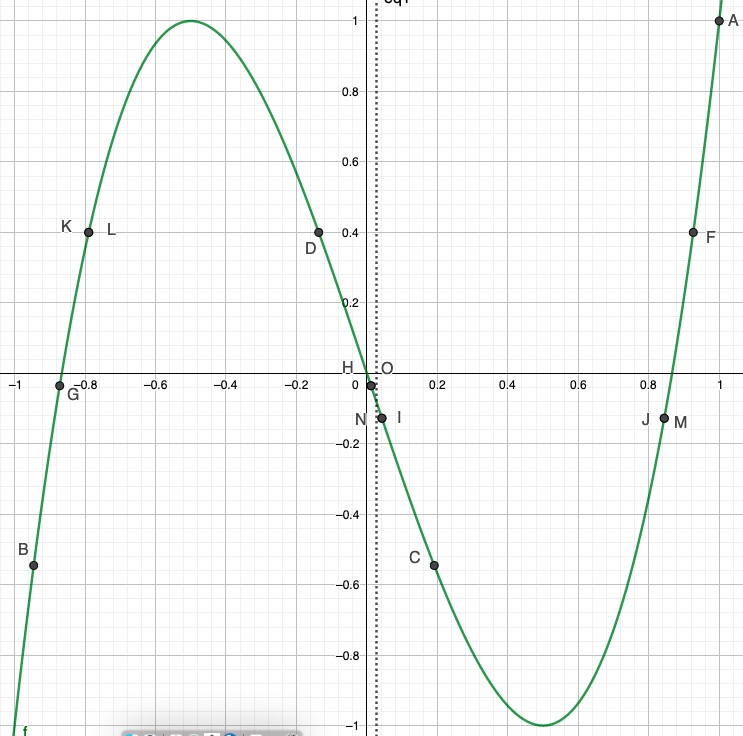}
                           \includegraphics[scale=0.15]{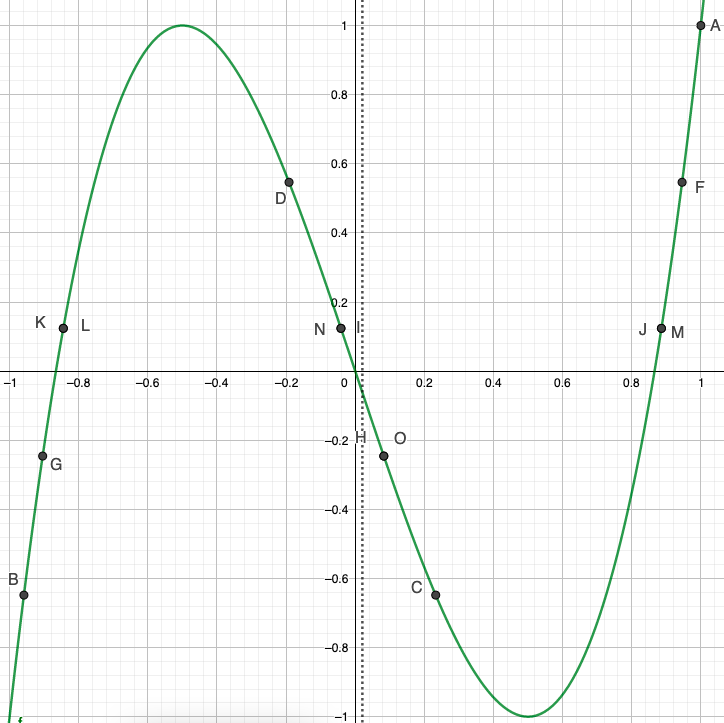}
                             \includegraphics[scale=0.15]{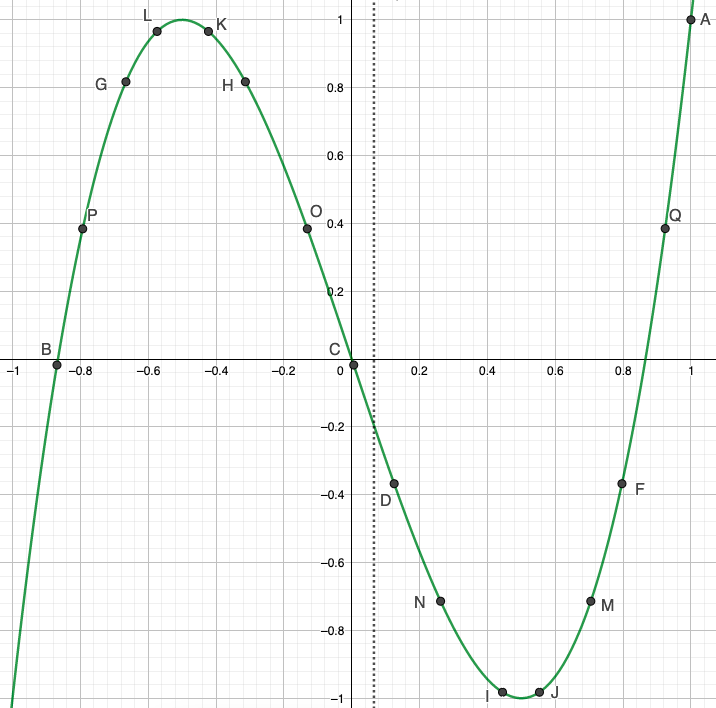}
                              \includegraphics[scale=0.15]{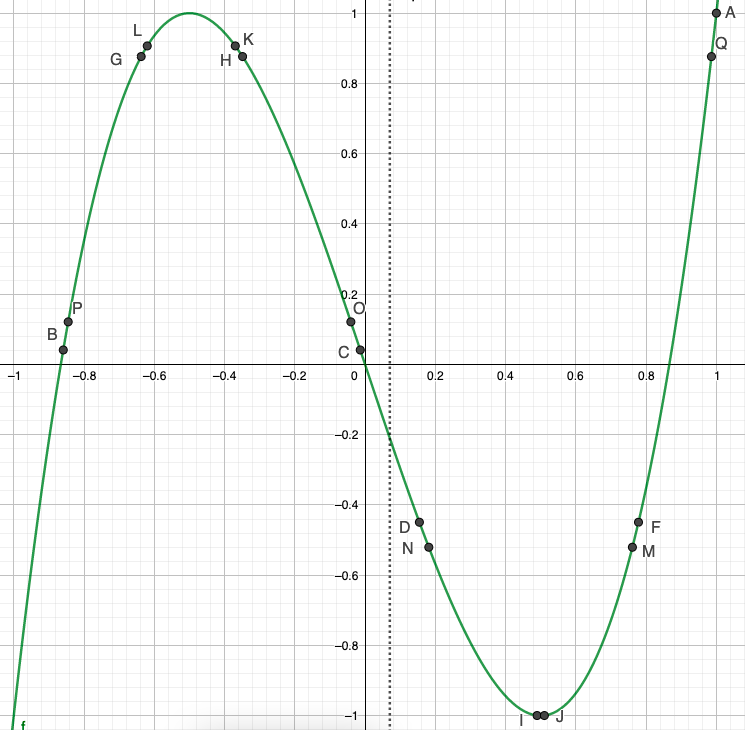}
                              \includegraphics[scale=0.15]{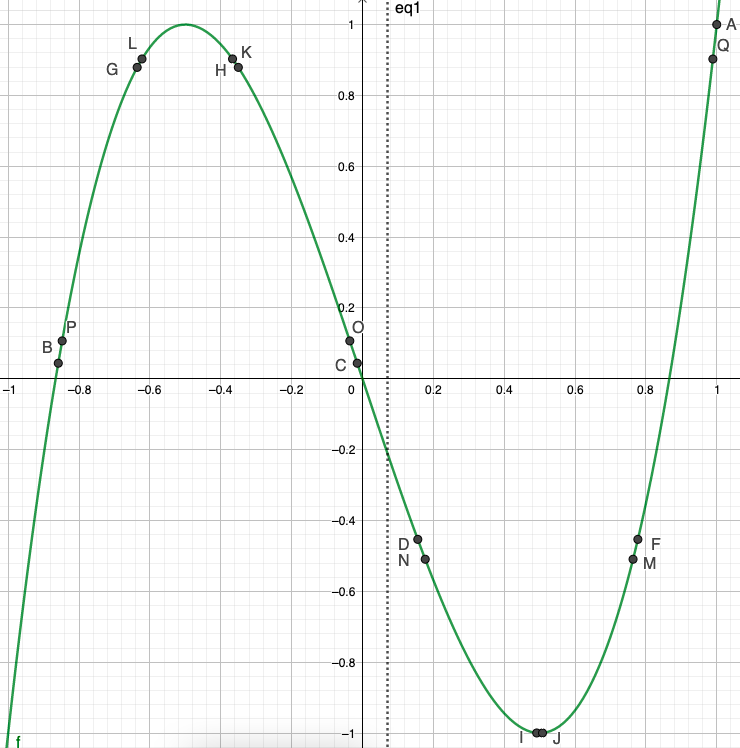}
                                \includegraphics[scale=0.15]{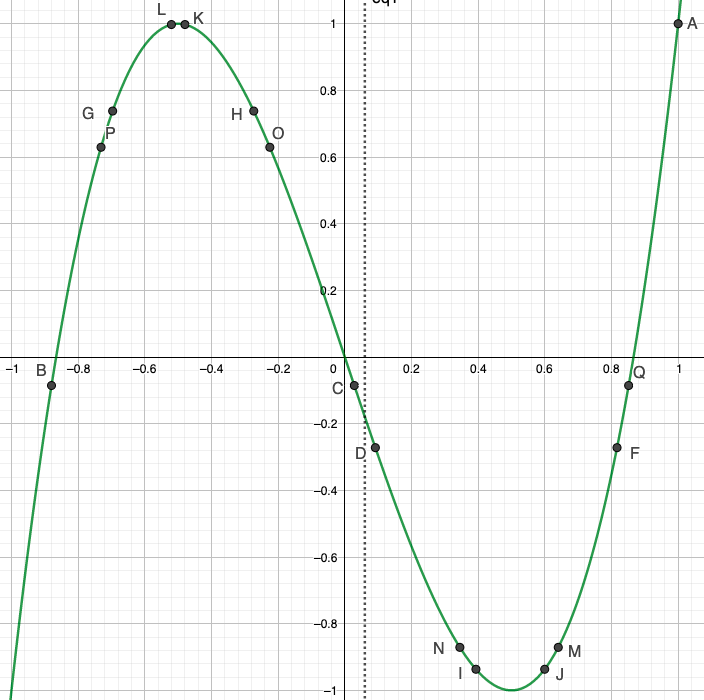}
                     \includegraphics[scale=0.15]{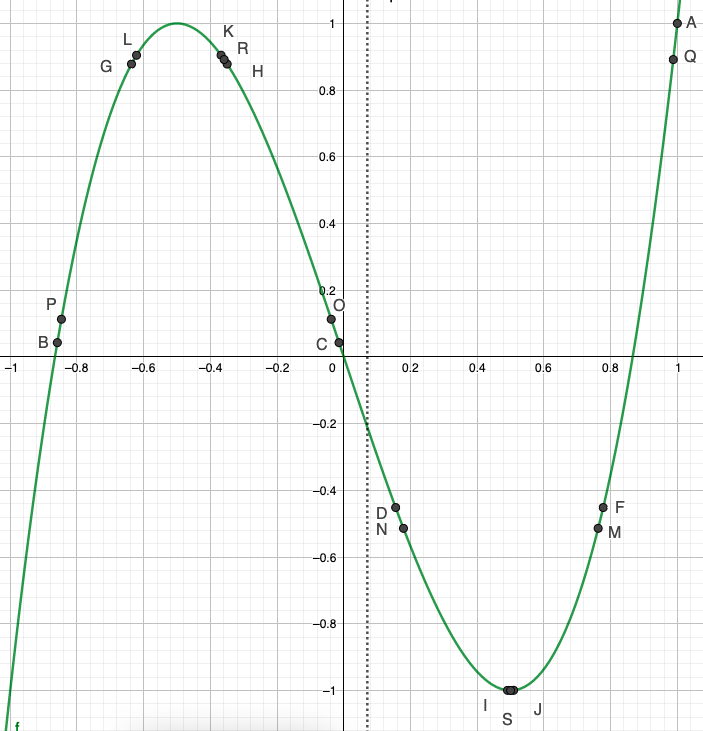}
         		\includegraphics[scale=0.15]{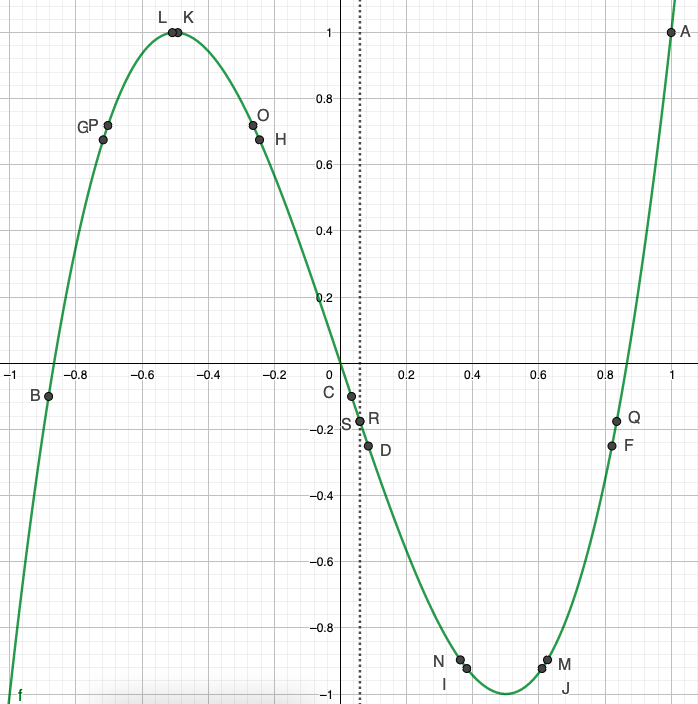}
		
\caption{}
\label{fig:work_in_progress2}
\end{figure}

\begin{figure}[H]
  \centering
  \includegraphics[scale=0.15]{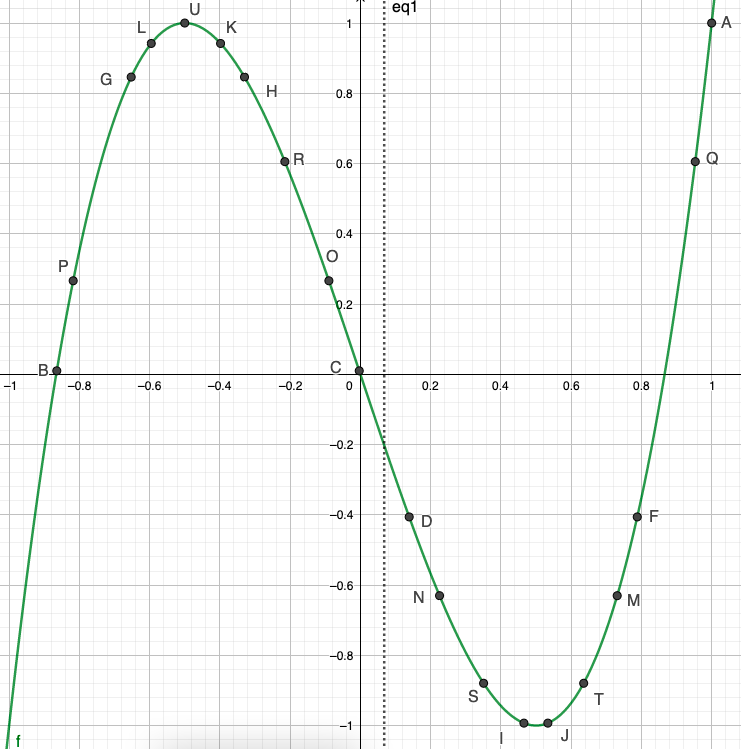}
   \includegraphics[scale=0.15]{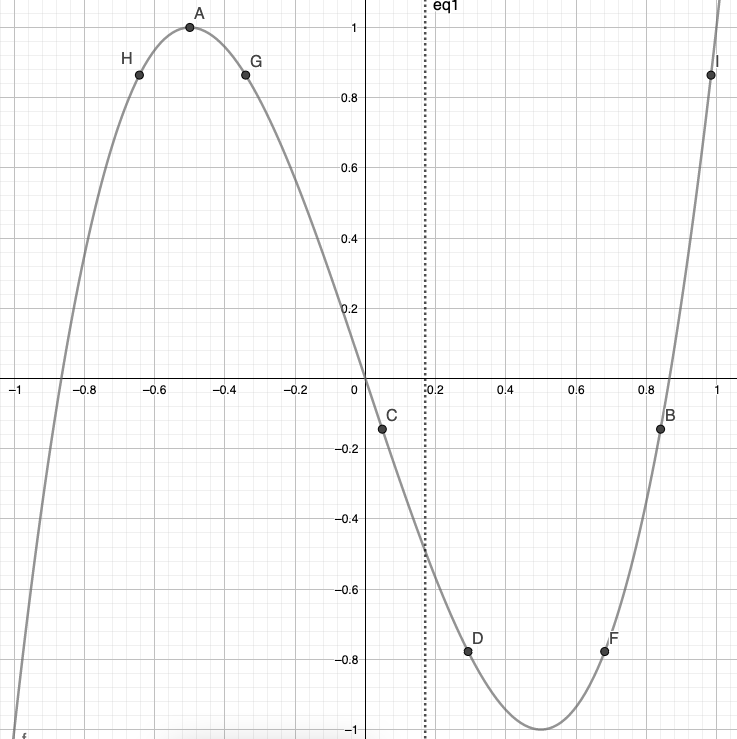}
      \includegraphics[scale=0.15]{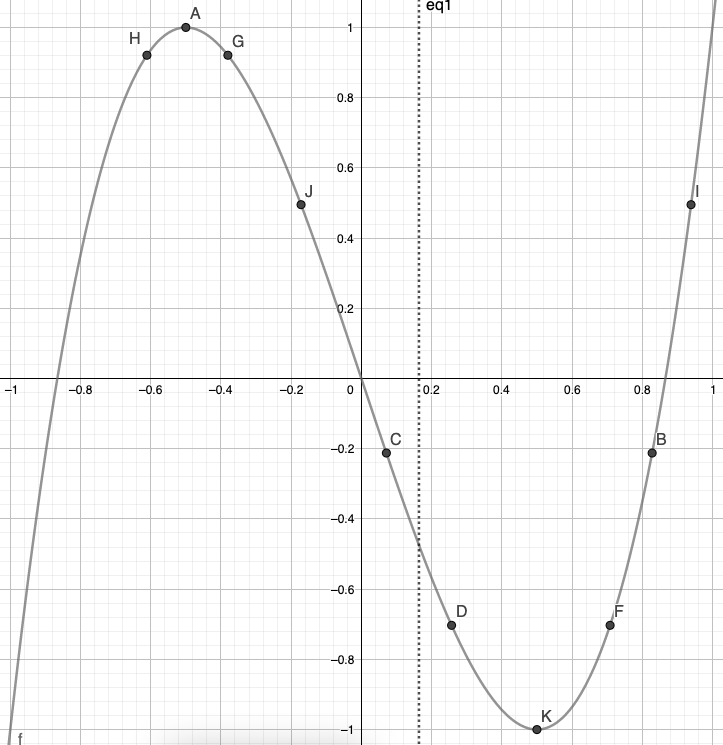}
      \includegraphics[scale=0.15]{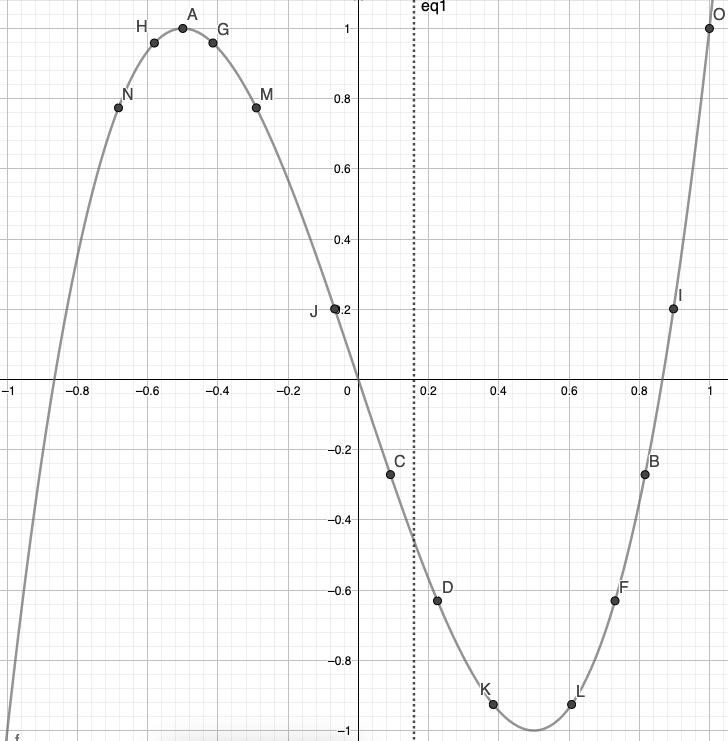}
     \includegraphics[scale=0.15]{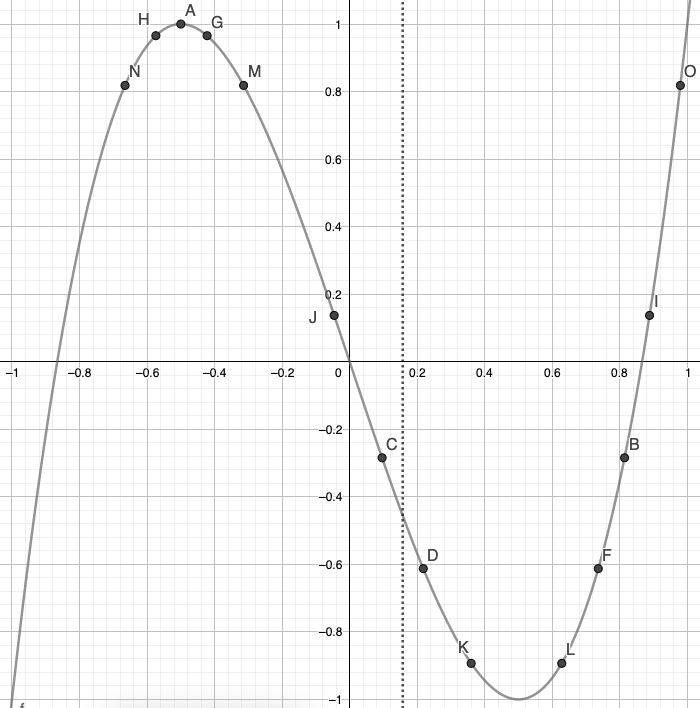}
         \includegraphics[scale=0.15]{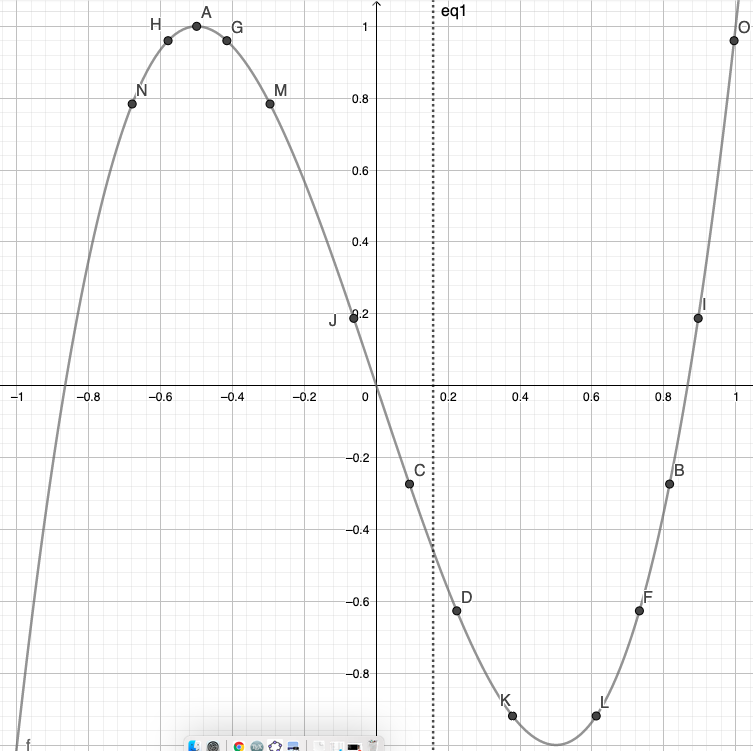}
          \includegraphics[scale=0.15]{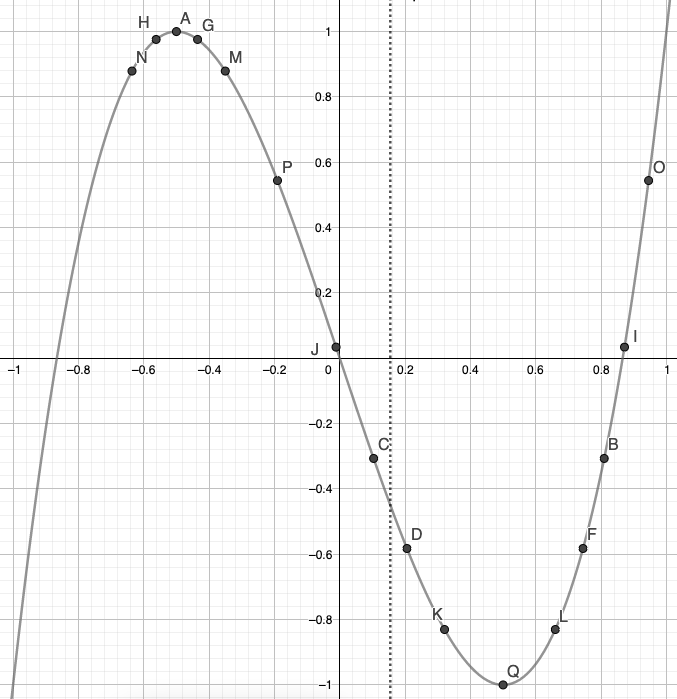}
            \includegraphics[scale=0.15]{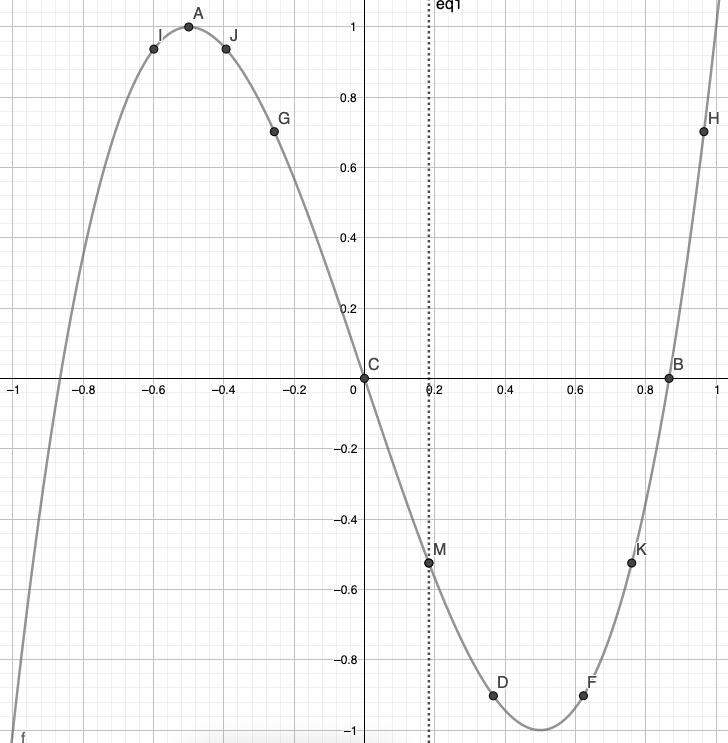}
            \includegraphics[scale=0.15]{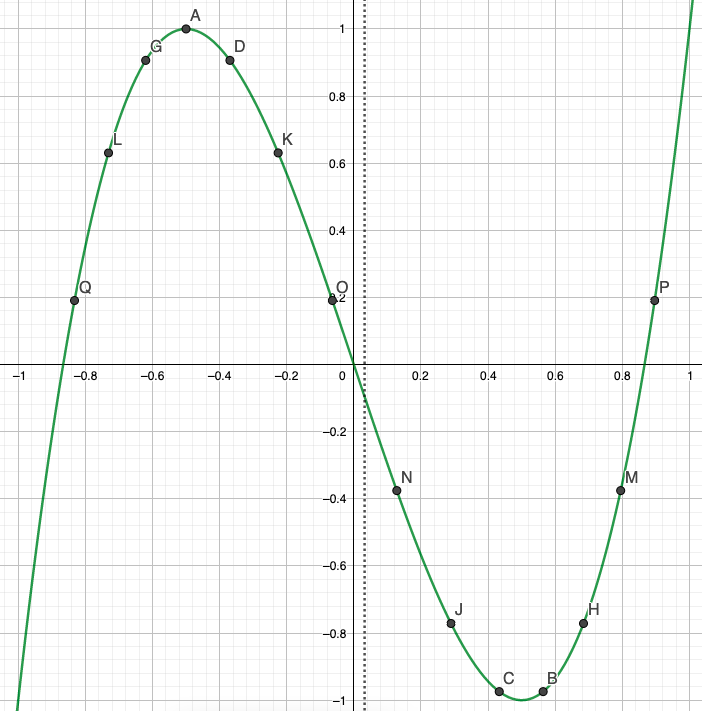}
             \includegraphics[scale=0.15]{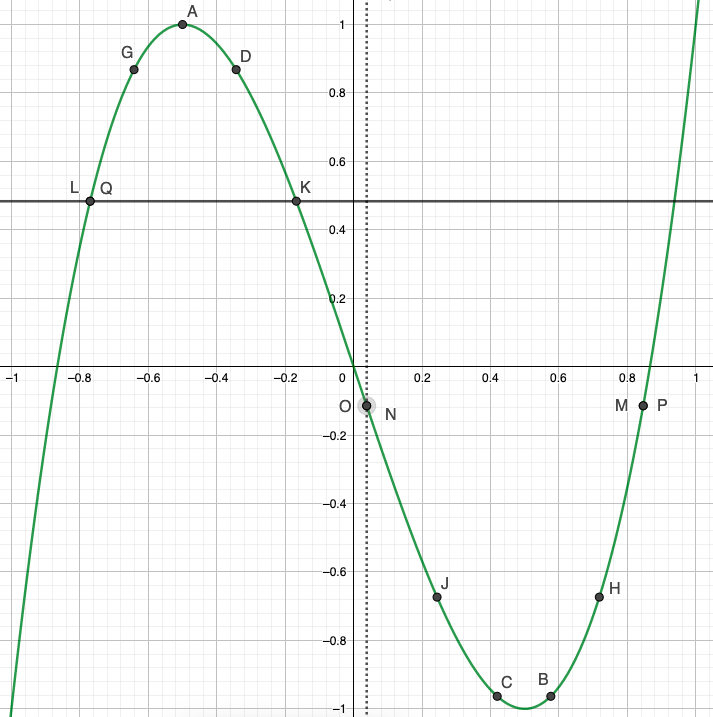}
              \includegraphics[scale=0.15]{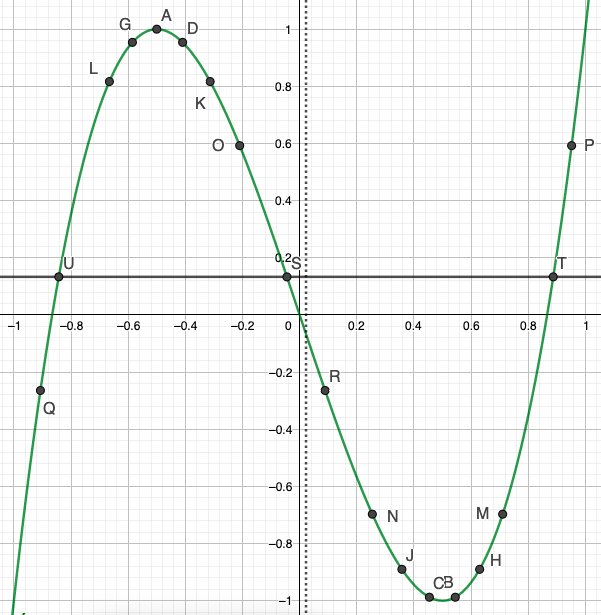}
               \includegraphics[scale=0.15]{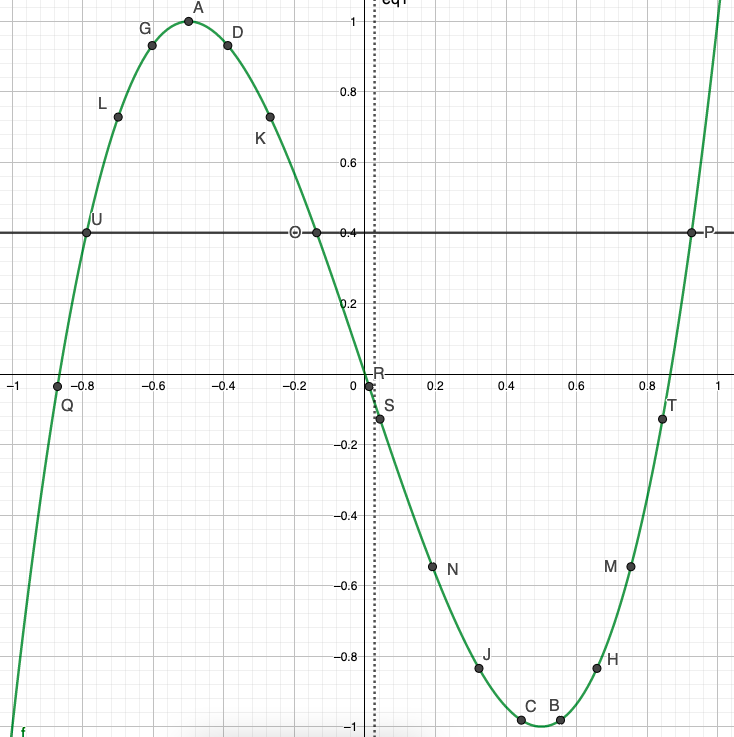}
                \includegraphics[scale=0.15]{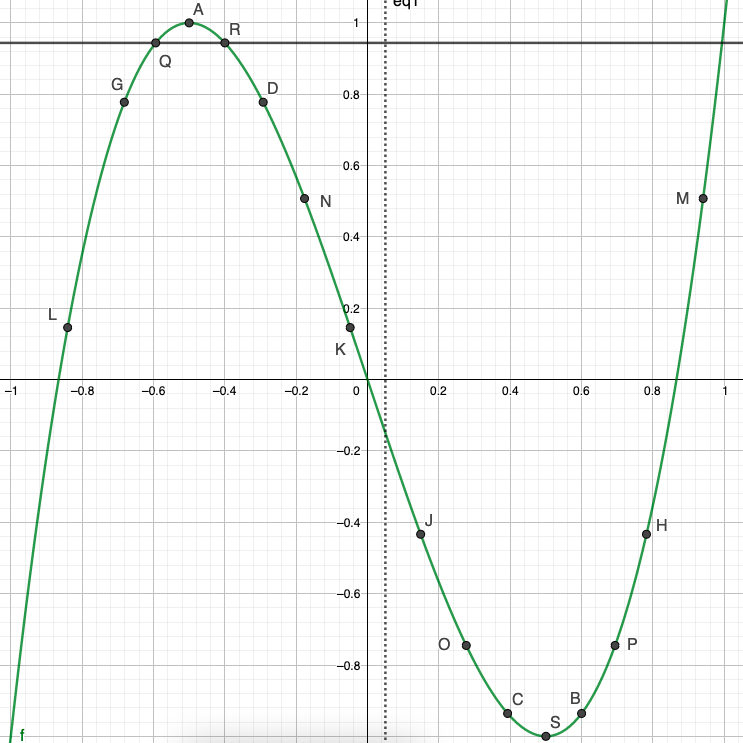}
                 \includegraphics[scale=0.15]{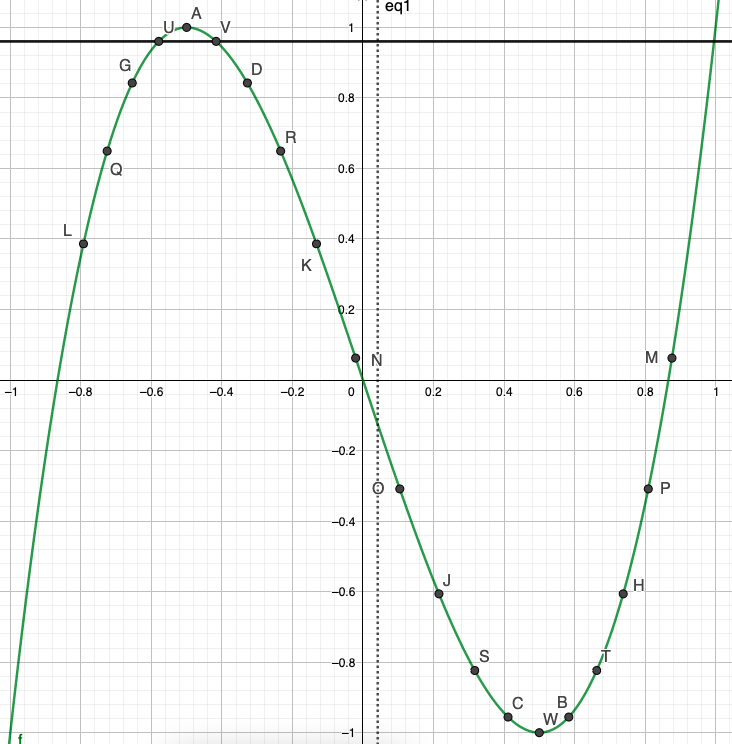}
            \includegraphics[scale=0.15]{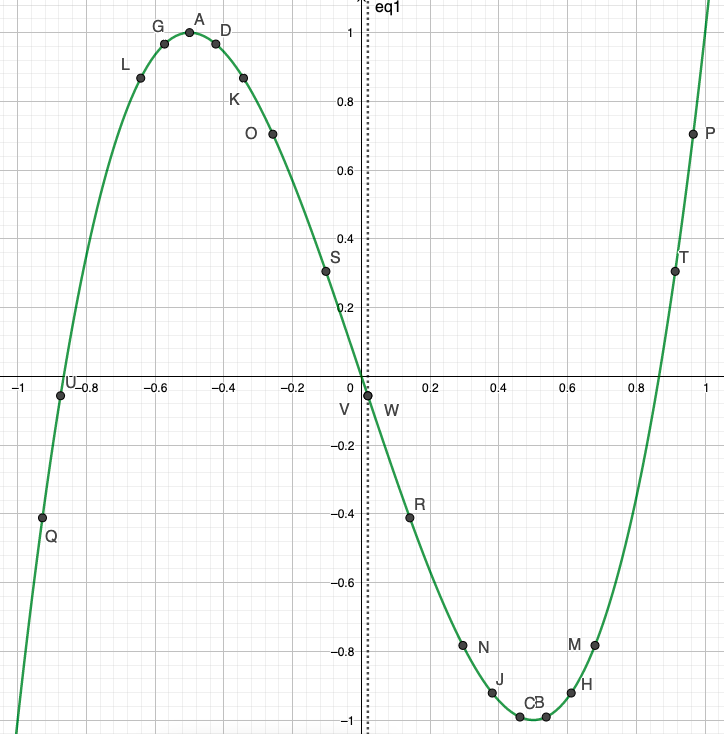}
            \includegraphics[scale=0.15]{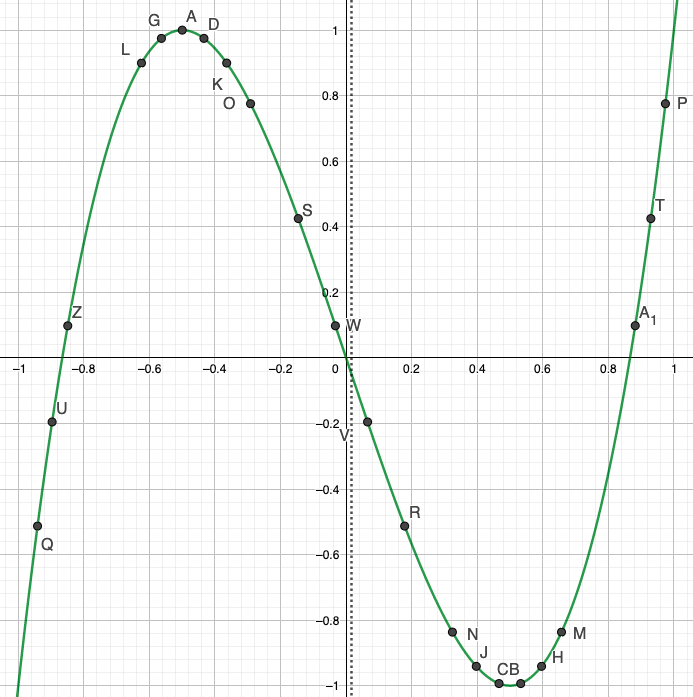}
            \includegraphics[scale=0.15]{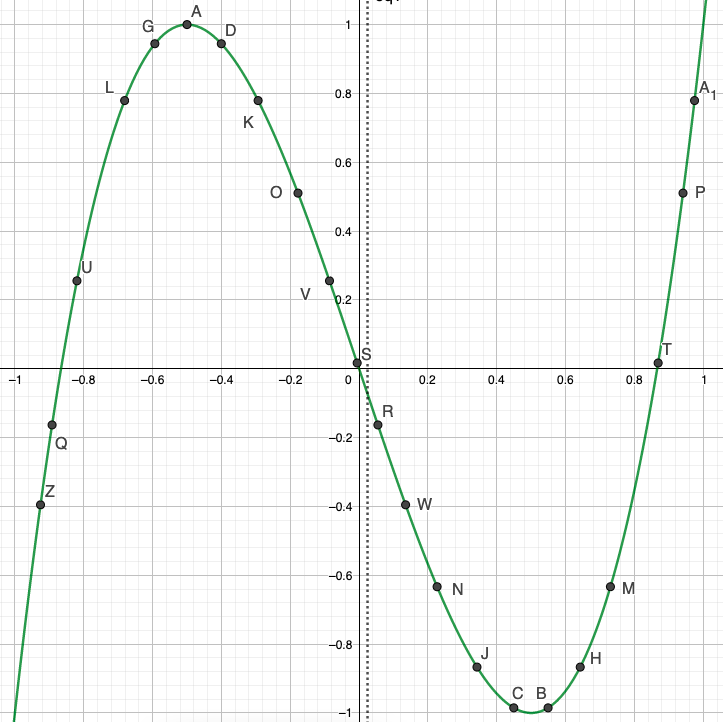}
             \includegraphics[scale=0.15]{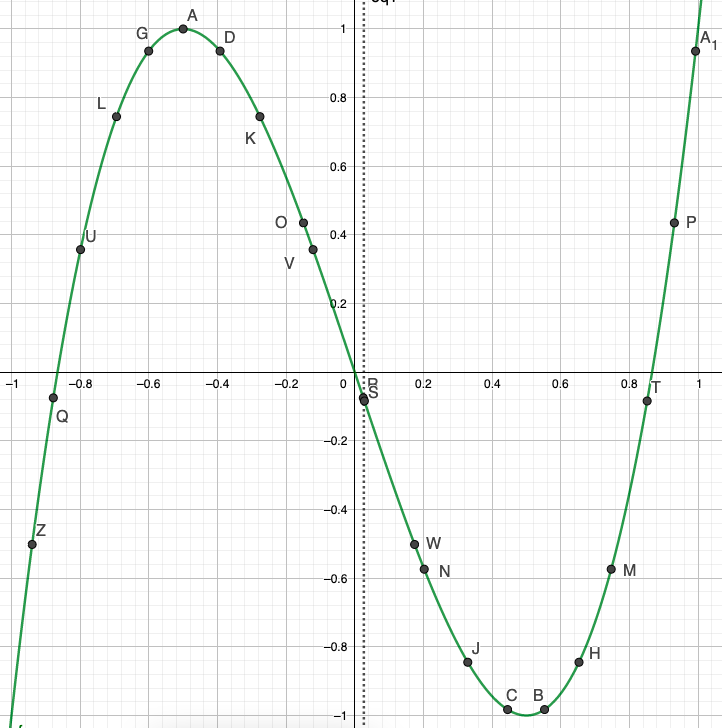}
              \includegraphics[scale=0.15]{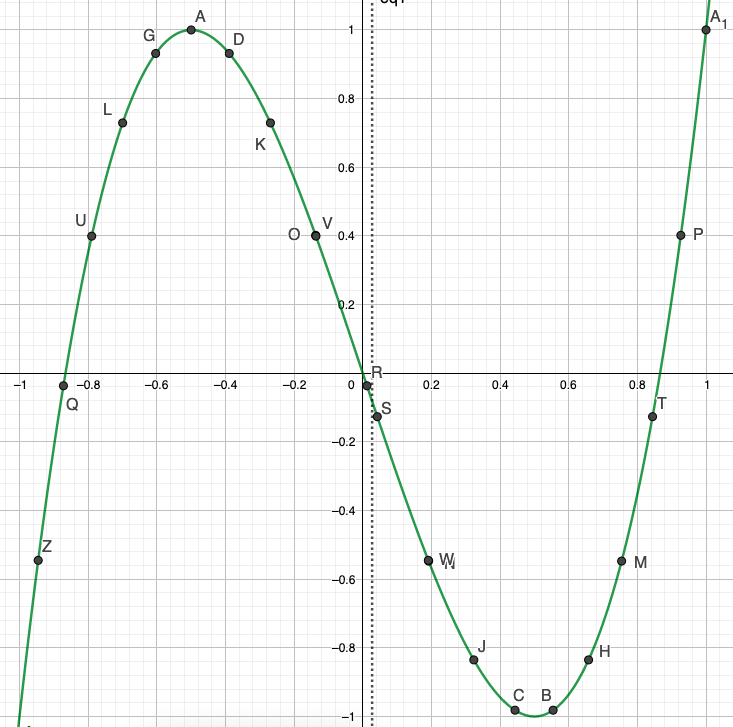}
               \includegraphics[scale=0.15]{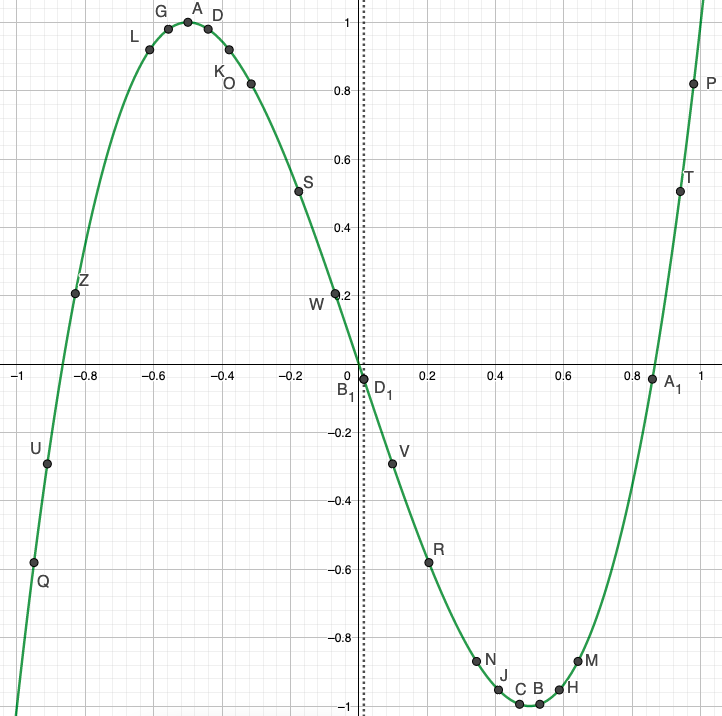}
                \includegraphics[scale=0.15]{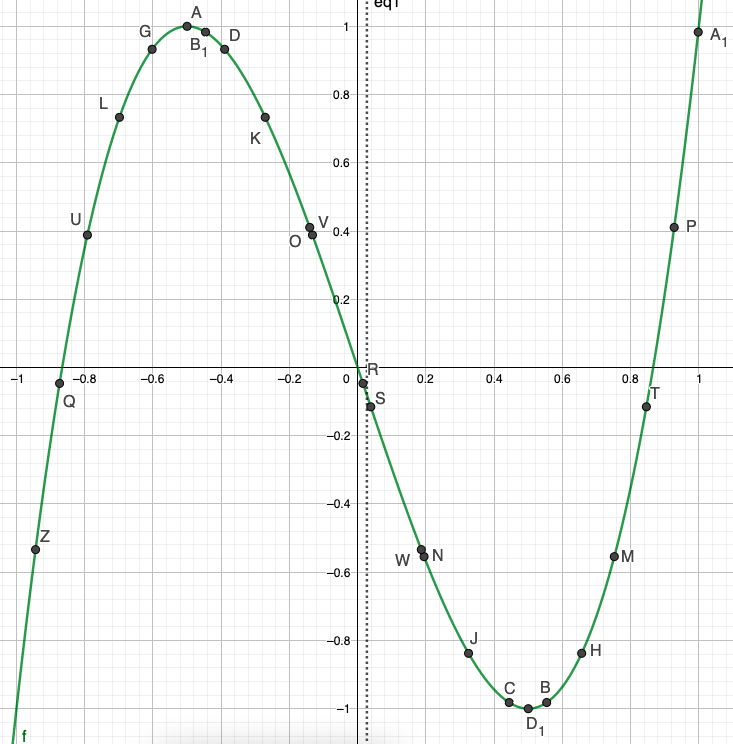}
                 \includegraphics[scale=0.15]{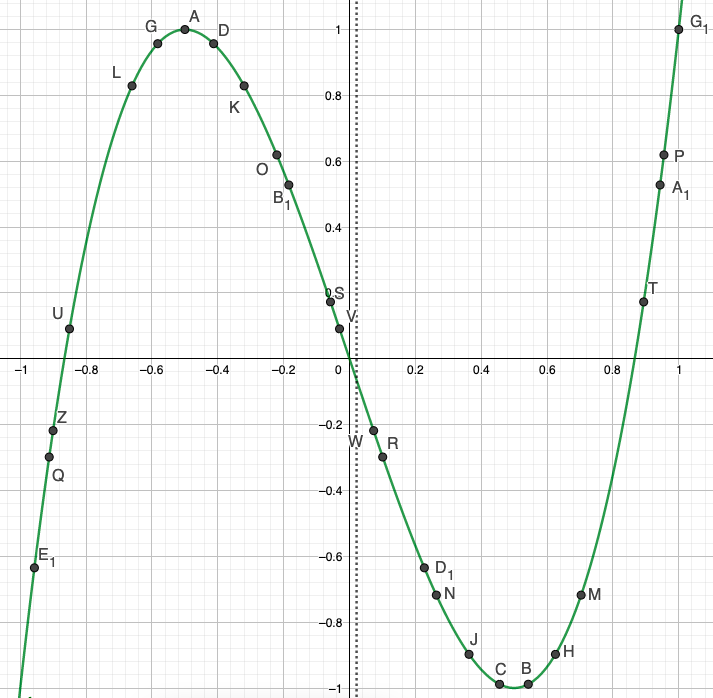}
                  \includegraphics[scale=0.15]{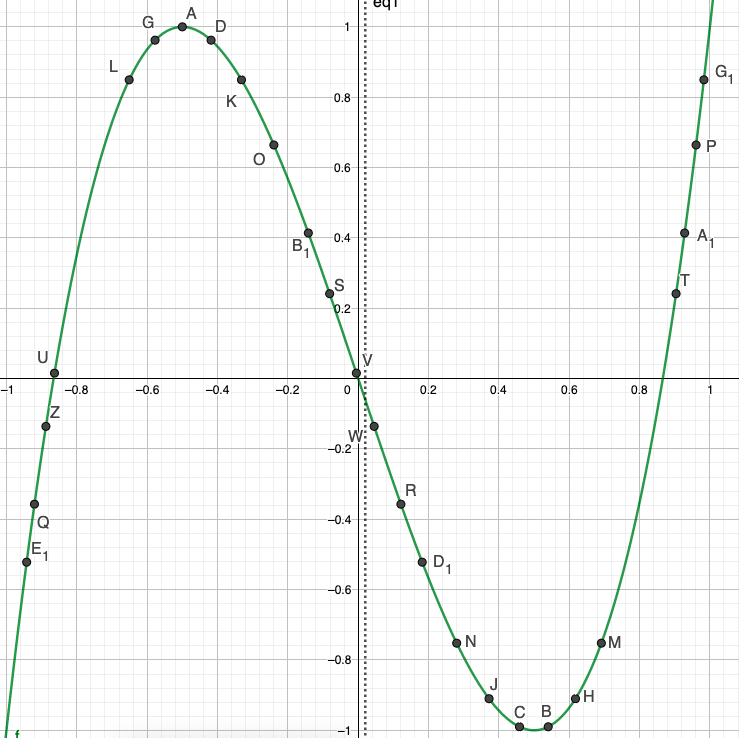}
                 
\caption{}
\label{fig:work_in_progress3}
\end{figure}

1) $\E \simeq 0.169$. 

2) $\E = 0.25$.

3) $\E = -0.25$.

4) $\E \simeq -0.1028$.

5) $\E \simeq -0.16019$.

6) $\E \simeq -0.1202$.

7) $\E \simeq -0.17911$.

8) $\E \simeq -0.11608$.

9) $\E \simeq -0.10796$.

10) $\E \simeq -0.11385$.

11) $\E \simeq -0.08425$.

12) $\E \simeq -0.0917$.

13) $\E \simeq -0.07143$.

14) $\E \simeq -0.14325$.

15) $\E \simeq -0.05169$.

16) $\E \simeq -0.06615$.

17) $\E \simeq -0.08424$.

18) $\E \simeq -0.18726$.

19) $\E \simeq -0.2857$.

20) $\E \simeq -0.40824$.

21) $\E \simeq -0.2857$.

22) $\E \simeq -0.13911$.

23) $\E \simeq -0.06805$.

24) $\E \simeq -0.0813$.

25) $\E \simeq -0.07704$.

26) $\E \simeq -0.05406$.

27) $\E \simeq -0.08014$.

28) $\E \simeq -0.04594$.

29) $\E \simeq -0.27129$.

30) $\E \simeq -0.1082$.

31) $\E \simeq -0.38229$.

32) $\E \simeq -0.4467$.

33) $\E \simeq -0.34518$.

34) $\E \simeq -0.43085$.

35) $\E \simeq -0.105384$.

36) $\E \simeq -0.05659$.

37) $\E \simeq -0.04657$.

38) $\E \simeq -0.09839$.

39) $\E \simeq -0.08979$.

40) $\E \simeq -0.03862$.

41) $\E \simeq +0.05439$.

42) $\E \simeq +0.04135$.

43) $\E \simeq +0.13104$.

44) $\E \simeq +0.14092$.

45) $\E \simeq +0.14142$.

46) $\E \simeq +0.12009$.

47) $\E \simeq +0.14121$.

48) $\E \simeq +0.11771$.

49) $\E \simeq +0.13562$.

50) $\E \simeq +0.34089$.

51) $\E \simeq +0.32821$.

52) $\E \simeq +0.31674$. 

53) $\E \simeq +0.31404$. 

54) $\E \simeq +0.31616$.

55) $\E \simeq +0.30948$.

56) $\E \simeq +0.36612$.

57) $\E \simeq +0.06397$.

58) $\E \simeq +0.07628$.

59) $\E \simeq +0.04454$.

60) $\E \simeq +0.05450$.

61) $\E \simeq +0.1001$.

62) $\E \simeq +0.08365$.

63) $\E \simeq +0.03773$.

64) $\E \simeq +0.03277$.

65) $\E \simeq +0.04901$.

66) $\E \simeq +0.05283$.

67) $\E \simeq +0.05445$.

68) $\E \simeq +0.02930$.

69) $\E \simeq +0.05402$.

70) $\E \simeq +0.04300$.

71) $\E \simeq +0.04030$.

\end{document}